\edef\inewcount{\noexpand\csname newcount\endcsname}
\edef\inewdimen{\noexpand\csname newdimen\endcsname}
\edef\inewskip{\noexpand\csname newskip\endcsname}
\edef\inewmuskip{\noexpand\csname newmuskip\endcsname}
\edef\inewbox{\noexpand\csname newbox\endcsname}
\edef\inewhelp{\noexpand\csname newhelp\endcsname}
\edef\inewtoks{\noexpand\csname newtoks\endcsname}
\edef\inewread{\noexpand\csname newread\endcsname}
\edef\inewwrite{\noexpand\csname newwrite\endcsname}
\edef\inewfam{\noexpand\csname newfam\endcsname}
\edef\inewlanguage{\noexpand\csname newlanguage\endcsname}
\edef\inewinsert{\noexpand\csname newinsert\endcsname}
\edef\inewif{\noexpand\csname newif\endcsname}

\countdef\ch=253
\ch="80 \loop\ifnum\ch<"100 \lccode\ch=0 \uccode\ch=0 \advance\ch1 \repeat % for LaTeX
\ch="80 \loop\ifnum\ch<"C0 \catcode\ch=11 \advance\ch1 \repeat % make all intermediate UTF-8 octets (top 2 bits are 10) letters
\ch="C0 \loop\ifnum\ch<"100 \catcode\ch=15 \advance\ch1 \repeat % make all octets with the top 2 bits set illegal
\ch="C2 \loop\ifnum\ch<"F5 \catcode\ch=\active \advance\ch1 \repeat % make all legitimate initial UTF-8 octets active
\ch="C2 \loop\ifnum\ch<"E2 \uccode\ch=1 \advance\ch1 \repeat % make uccode=1 for all initial UTF-8 octets of letters, to allow macros to recognize Unicode letters
\let\ch\underfined

\catcode0=12 \def\cdef#1#2{\begingroup\lccode0=`#2 \lowercase{\endgroup \def#1{^^@}}} \catcode0=9

\catcode`@=\active
\def@#1{\cdef\chr#1 \edef#1##1{\noexpand\csname\chr##1\noexpand\endcsname}} % 2-byte
@^^c2@^^c3@^^c4@^^c5@^^c6@^^c7@^^c8@^^c9@^^ca@^^cb@^^cc@^^cd@^^ce@^^cf@^^d0@^^d1@^^d2@^^d3@^^d4@^^d5@^^d6@^^d7@^^d8@^^d9@^^da@^^db@^^dc@^^dd@^^de@^^df
\def@#1{\cdef\chr#1 \edef#1##1##2{\noexpand\csname\chr##1##2\noexpand\endcsname}} % 3-byte
@^^e0@^^e1@^^e2@^^e3@^^e4@^^e5@^^e6@^^e7@^^e8@^^e9@^^ea@^^eb@^^ec@^^ed@^^ee@^^ef
\def@#1{\cdef\chr#1 \edef#1##1##2##3{\noexpand\csname\chr##1##2##3\noexpand\endcsname}} % 4-byte
@^^f0@^^f1@^^f2@^^f3@^^f4
\let\chr\undefined
\let\cdef\undefined

\def\grabfuturelet{\futurelet\next\grabexamine}
\def\grabexamine{\ifx\next\csname\expandafter\grab\fi}
\obeylines \def\grab\csname#1\endcsname#2^^M{\expandafter\def\csname#1\endcsname{#2}\expandafter\grabfuturelet} \expandafter\grabfuturelet%
 ~
¢{\hbox{\rm\rlap/c}}
£{\it\$}
«\leftguillemet
­\-
»\rightguillemet
À{\`A}
Á{\'A}
Â{\^A}
Ã{\~A}
Ä{\"A}
Ç{\c C}
È{\`E}
É{\'E}
Ê{\^E}
Ë{\"E}
Ì{\`I}
Í{\'I}
Î{\^I}
Ï{\"I}
Ð\ETH
Ñ{\~N}
Ò{\`O}
Ó{\'O}
Ô{\^O}
Õ{\~O}
Ö{\"O}
Ù{\`U}
Ú{\'U}
Û{\^U}
Ü{\"U}
Ý{\'Y}
Þ\THORN
à{\`a}
á{\'a}
â{\^a}
ã{\~a}
ä{\"a}
ç{\c c}
è{\`e}
é{\'e}
ê{\^e}
ë{\"e}
ì{\`\i}
í{\'\i}
î{\^\i}
ï{\"\i}
ð\eth
ñ{\~n}
ò{\`o}
ó{\'o}
ô{\^o}
õ{\~o}
ö{\"o}
ù{\`u}
ú{\'u}
û{\^u}
ü{\"u}
ý{\'y}
þ\thorn
ÿ{\"y}
Ā{\=A}
ā{\=a}
Ă{\u A}
ă{\u a}
Ć{\'C}
ć{\'c}
Ĉ{\^C}
ĉ{\^c}
Ċ{\.C}
ċ{\.c}
Č{\v C}
č{\v c}
Ď{\v D}
ď{\v d}
Ē{\=E}
ē{\=e}
Ĕ{\u E}
ĕ{\u e}
Ė{\.E}
ė{\.e}
Ě{\v E}
ě{\v e}
Ĝ{\^G}
ĝ{\^g}
Ğ{\u G}
ğ{\u g}
Ġ{\.G}
ġ{\.g}
Ģ{\c G}
ģ{\c g}
Ĥ{\^H}
ĥ{\^h}
Ĩ{\~I}
ĩ{\~\i}
Ī{\=I}
ī{\=\i}
Ĭ{\u I}
ĭ{\u\i}
İ{\.I}
ĲIJ
ĳij
Ĵ{\^J}
ĵ{\^\j}
Ķ{\c K}
ķ{\c k}
Ĺ{\'L}
ĺ{\'l}
Ļ{\c L}
ļ{\c l}
Ľ{\v L}
ľ{\v l}
Ń{\'N}
ń{\'n}
Ņ{\c N}
ņ{\c n}
Ň{\v N}
ň{\v n}
Ō{\=O}
ō{\=o}
Ŏ{\u O}
ŏ{\u o}
Ő{\"O}
ő{\"o}
Ŕ{\'R}
ŕ{\'r}
Ŗ{\c R}
ŗ{\c r}
Ř{\v R}
ř{\v r}
Ś{\'S}
ś{\'s}
Ŝ{\^S}
ŝ{\^s}
Ş{\c S}
ş{\c s}
Š{\v S}
š{\v s}
Ţ{\c T}
ţ{\c t}
Ť{\v T}
ť{\v t}
Ũ{\~U}
ũ{\~u}
Ū{\=U}
ū{\=u}
Ŭ{\u U}
ŭ{\u u}
Ű{\H U}
ű{\H u}
Ŵ{\^W}
ŵ{\^w}
Ŷ{\^Y}
ŷ{\^y}
Ÿ{\"Y}
Ź{\'Z}
ź{\'z}
Ż{\.Z}
ż{\.z}
Ž{\v Z}
ž{\v z}
’'
‘`
”{''}
“{``}
‐-
–{--}
—{---}
¡{!`}
¿{?`}
−-
′'
ß\ss
æ\ae
Æ\AE
œ\oe
Œ\OE
ø\o
Ø\O
å\aa
Å\AA
ł\l
Ł\L
†\dag
‡\ddag
§\S
¶\P
©\copyright
…\dots
ı\relax\ifmmode\imath\else\i\fi
ȷ\relax\ifmmode\jmath\else\j\fi
¬\neg
△\bigtriangleup
÷\div
±\pm
◯\bigcirc
×\times
‖\|
←\leftarrow
→\rightarrow
∥\Vert
⟩\rangle
⟨\langle
∣|
∮\oint
≐\dot=
↑\uparrow
↓\downarrow
α\alpha
β\beta
γ\gamma
δ\delta
ϵ\epsilon
ζ\zeta
η\eta
θ\theta
ι\iota
κ\kappa
λ\lambda
μ\mu
ν\nu
ξ\xi
οo
π\pi
ρ\rho
σ\sigma
τ\tau
υ\upsilon
ϕ\phi
χ\chi
ψ\psi
ω\omega
ε\varepsilon
ϑ\vartheta
ϖ\varpi
ϱ\varrho
ς\varsigma
φ\varphi
Γ\Gamma
Δ\Delta
Θ\Theta
Λ\Lambda
Ξ\Xi
Π\Pi
Σ\Sigma
Υ\Upsilon
Φ\Phi
Ψ\Psi
Ω\Omega
ℵ\aleph
ℏ\hbar
ℓ\ell
℘\wp
ℜ\Re
ℑ\Im
∂\partial
∞\infty
∅\emptyset
∇\nabla
√\surd
⊤\top
⊥\bot
∠\angle
△\triangle
∀\forall
∃\exists
♭\flat
♮\natural
♯\sharp
♣\clubsuit
♢\diamondsuit
♡\heartsuit
♠\spadesuit
∐\coprod
⋁\bigvee
⋀\bigwedge
⨄\biguplus
⋂\bigcap
⋃\bigcup
∫\int
∏\prod
∑\sum
⨂\bigotimes
⨁\bigoplus
⨀\bigodot
⨆\bigsqcup
◁\triangleleft
▷\triangleright
▽\bigtriangledown
∧\wedge
∨\vee
∩\cap
∪\cup
⊓\sqcap
⊔\sqcup
⊎\uplus
⨿\amalg
⋄\diamond
∙\bullet
≀\wr
⊙\odot
⊘\oslash
⊗\otimes
⊖\ominus
⊕\oplus
∓\mp
∘\circ
○\Orb
∖\setminus
⋅\cdot
∗\ast
⨯\times
⋆\star
∝\propto
⊑\sqsubseteq
⊒\sqsupseteq
∥\parallel
∣\divides
⊣\dashv
⊢\vdash
↗\nearrow
↘\searrow
↖\nwarrow
↙\swarrow
⇔\Leftrightarrow
⇐\Leftarrow
⇒\Rightarrow
≠\neq
≤\leq
≥\geq
≻\succ
≺\prec
≈\approx
≽\succeq
≼\preceq
⊃\supset
⊂\subset
⊇\supseteq
⊆\subseteq
∈\in
∋\ni
≫\gg
≪\ll
↔\leftrightarrow
↦\mapsto
∼\sim
≃\simeq
⟂\perp
≡\equiv
≍\asymp
⌣\smile
⌢\frown
↼\leftharpoonup
↽\leftharpoondown
⇀\rightharpoonup
⇁\rightharpoondown
↪\hookrightarrow
↩\hookleftarrow
⋈\bowtie
⊨\models
⟹\Longrightarrow
⟶\longrightarrow
⟵\longleftarrow
⟸\Longleftarrow
⟼\longmapsto
⟷\longleftrightarrow
⟺\Longleftrightarrow
⋯\cdots
⋮\vdots
⋱\ddots
↕\updownarrow
⇑\Uparrow
⇓\Downarrow
⇕\Updownarrow
⌉\rceil
⌈\lceil
⌋\rfloor
⌊\lfloor
≅\cong
∉\notin
⇌\rightleftharpoons
≐\doteq
∄\not\exists
∌\not\ni
∔\dot+
∕/
∤\not|
∦\not\|
∬\int\!\!\!\int
∭\int\!\!\!\int\!\!\!\int
∮\oint
∸\dot-
≁\not\sim
≄\not\simeq
≆\not\cong
≇\not\cong
≉\not\approx
≔:=
≕=:
≢\not\equiv
≭\not\asump
≮\not<
≯\not<
≰\not\le
≱\not\ge
⊀\not\prec
⊁\not\succ
⊄\not\subset
⊅\not\supset
⊈\not\subseteq
⊉\not\supseteq
⊦\vdash
⊧\models
⊬\not\vdash
⊭\not\models
⊲\triangleleft
⊳\triangleright
⋠\not\preceq
⋡\not\succeq
⋤\not\sqsubseteq
⋥\not\sqsupseteq
⋪\not\triangleleft
⋫\not\triangleright
◻\square

\catcode`\^^M=5 %
\let\grabfuturelet\undefined \let\grabexamine\undefined \let\grab\undefined

\let\xcsname=\csname
\let\xendcsname=\endcsname
\def@#1{\def#1##1{\expandafter\ifx\csname\string#1##1\endcsname\relax\errmessage{Undefined UTF-8 sequence \string#1##1}\else\xcsname\string#1##1\xendcsname\fi}}
@^^c2@^^c3@^^c4@^^c5@^^c6@^^c7@^^c8@^^c9@^^ca@^^cb@^^cc@^^cd@^^ce@^^cf@^^d0@^^d1@^^d2@^^d3@^^d4@^^d5@^^d6@^^d7@^^d8@^^d9@^^da@^^db@^^dc@^^dd@^^de@^^df
\def@#1{\def#1##1##2{\expandafter\ifx\csname\string#1##1##2\endcsname\relax\errmessage{Undefined UTF-8 sequence \string#1##1##2}\else\xcsname\string#1##1##2\xendcsname\fi}}
@^^e0@^^e1@^^e2@^^e3@^^e4@^^e5@^^e6@^^e7@^^e8@^^e9@^^ea@^^eb@^^ec@^^ed@^^ee@^^ef
\def@#1{\def#1##1##2##3{\expandafter\ifx\csname\string#1##1##2##3\endcsname\relax\errmessage{Undefined UTF-8 sequence \string#1##1##2##3}\else\xcsname\string#1##1##2##3\xendcsname\fi}}
@^^f0@^^f1@^^f2@^^f3@^^f4
\let@\undefined
\catcode`@=12

\newif\ifscroll % scroll vs codex
\newif\ifsuppressunusedbib % suppress unused bibliography items?

\def\printerr#1{\immediate\write17{#1}}
\def\warningline#1#2{\printerr{! #2}\printerr{l.#1}\printerr{}}
\def\ewarningline#1#2#3{\printerr{! #2}\printerr{l.#1 #3}\printerr{}}
\def\warning{\warningline{\the\inputlineno}}
\long\def\gobble#1{}
\ifx\gobbleinit\undefined{\long\gdef\gobbleinit#1\par{}}\fi
\def\expand#1{\edef\expandmacro{#1}\expandmacro\let\expandmacro\undefined}
\def\setetok#1#2{\expand{\noexpand#1{#2}}}
\def\expandtoks#1{\expandafter\edef\expandafter\expandmacro\expandafter{\the#1}#1\expandafter{\expandmacro}}
\def\appendexpand#1#2{\setetok#1{\the#1#2}}
\long\def\append#1#2{#1\expandafter{\the#1#2}}
\long\def\appendtoksexpand#1#2{#1\expandafter\expandafter\expandafter{\expandafter\the\expandafter#1\the#2}}
\long\def\appendonceexpand#1#2{#1\expandafter\expandafter\expandafter{\expandafter\the\expandafter#1#2}}

\def\link#1#2{\lhighlight{#2}}
\def\llink#1{\printlink{llink #1}\link{\ohash#1}}
\catcode`\#=11 \def\ohash{#}\catcode`\#=6
\catcode`\&=11 \def\ampersand{&}\catcode`\&=4
\def\anchor#1#2{\printlink{anchor #1 #2}#2}

\def\setpapersize#1#2{} % \number#1 sp triggers bugs in dvips and does not work with dvipdfm(x)
\def\dumpbox#1#2#3{\shipout\vbox{\setpapersize{#1}{#2}\unvbox#3}}
\def\mps#1{\epsfbox{#1}}
\def\metadata#1#2{}
\def\src{} % for proper match between bibliographic references and source code

\newread\epsffilein
\newif\ifepsfbbfound\inewif\ifepsffilecont
\newdimen\epsfxsize\inewdimen\epsfysize
\newdimen\pspoints\pspoints1bp
\let\runmp\errmessage % will be set to \warning if the METAPOST file cannot be compiled
\def\epsfbox#1{\openin\epsffilein=#1 \ifeof\epsffilein\runmp{Could not open file #1}\else
	{\def\do##1{\catcode`##1=12}\dospecials\catcode`\ =10\epsffileconttrue
		\epsfbbfoundfalse
		\loop\read\epsffilein to\epsffileline \ifeof\epsffilein\epsffilecontfalse\else\expandafter\epsfaux\epsffileline :. \\\fi\ifepsffilecont\repeat
		\ifepsfbbfound\else\errmessage{No HiResBoundingBox comment found in file #1}\fi}%
	\closein\epsffilein
	\epsfysize\epsfury\pspoints \advance\epsfysize-\epsflly\pspoints
	\epsfxsize\epsfurx\pspoints \advance\epsfxsize-\epsfllx\pspoints
	\setbox0\hbox{\vbox to\epsfury\pspoints{\vfil\hbox to\epsfxsize{\dimen0=\epsfllx\pspoints \kern-\dimen0 \includegraphics{#1}\hfil}}}%
	\dp0=\epsflly\pspoints \dp0=-\dp0
	\box0 \fi}
\catcode`\%=12 \def\epsfbblit{%%HiResBoundingBox} \catcode`\%=14
\let\dummybrace=} % to compensate for { on the previous line when it is read by \read
\def\epsfaux#1:#2\\{\def\testit{#1}\ifx\testit\epsfbblit \epsfgrab #2 . . . \\\epsffilecontfalse\epsfbbfoundtrue\fi}
\def\empty{}
\def\epsfgrab #1 #2 #3 #4 #5\\{\gdef\epsfllx{#1}\ifx\epsfllx\empty\epsfgrab #2 #3 #4 #5 .\\\else\gdef\epsflly{#2}\gdef\epsfurx{#3}\gdef\epsfury{#4}\fi} % gdef because epxfbox wraps us in a group

\newif\ifpdf \pdffalse \ifx\pdfoutput\undefined\else\ifx\pdfoutput\relax\else\ifnum\pdfoutput<1 \else\pdftrue\fi\fi\fi
\ifpdf
\def\pdflink#1#2{\leavevmode \lhighlight{\pdfstartlink user { /Subtype /Link /Border [0 0 0] /A << /S #1 >> }#2\pdfendlink}}
\def\llink#1{\pdflink{/GoTo /D (#1)}}
\def\link#1{\pdflink{/URI /URI (#1)}}
\def\anchor#1#2{\pdfdest name {#1} xyz #2}

\def\setpapersize#1#2{\pdfpagewidth#1 \pdfpageheight#2 }
\def\dumpbox#1#2#3{\setpapersize{#1}{#2}\shipout\box#3}
\def\metadata#1#2{\pdfinfo{/Title (#1) /Author (#2)}}
\input supp-pdf.mkii % the only external dependence of DPMAC: a converter from METAPOST's output to PDF
\def\mps#1{\convertMPtoPDF{#1}{1}{1}}
\fi

\newtoks\buffertoks
\def\addcode{\immediate\write\mpout}
\def\addunexpandedcode#1{{\toks0={#1}\addcode{\the\toks0}}}
\def\addcodebuffer#1{\edef\tmp{#1}\buffertoks\expandafter\expandafter\expandafter{\expandafter\the\expandafter\buffertoks\tmp}}
\def\addunexpandedcodebuffer#1{\buffertoks\expandafter{\the\buffertoks#1}}

\def\grabcode{\catcode`\#=12 \endlinechar=10 
	\afterassignment\dumpcode\outtoks} % make sure the newline characters get written
\def\dumpcode{\addcode{\the\outtoks}\endinlinemp\gobble} % \gobble the newline with character code 10
\def\begininlinemp{\inimp\begingroup\catcode`\^=7 \iftypesetting\mps{\filestem.\the\figno}\let\addcode\gobble\fi \addcode{beginfig(\the\figno);}}
\def\endinlinemp{\addcode{endfig;}\addcode{}\endgroup\global\advance\figno1\relax}
\ifx\endprolog\undefined\let\endprolog\relax\fi

\def\plainfmtname{plain}\ifx\fmtname\plainfmtname\else
\edef\plainoutput{\the\output}
\global\chardef\itfam=4
\def\_{\leavevmode \kern.06em \vbox{\hrule width.3em}} % for LaTeX

\font\tensy=cmsy10

\catcode"18=12
\catcode`@=11
{
\global\let\end\@@end
\global\let\input\@@input
}
\def\eqalign#1{\null\,\vcenter{\openup\jot\m@th
  \ialign{\strut\hfil$\displaystyle{##}$&$\displaystyle{{}##}$\hfil
      \crcr#1\crcr}}\,}
\catcode`@=12
\fi

\font\tenmsa=msam10 \font\sevenmsa=msam7 \font\fivemsa=msam5 \newfam\msafam \textfont\msafam=\tenmsa \scriptfont\msafam=\sevenmsa \scriptscriptfont\msafam=\fivemsa % 8
\font\teneufm=eufm10 \font\seveneufm=eufm7 \font\fiveeufm=eufm5 \newfam\eufmfam \textfont\eufmfam=\teneufm \scriptfont\eufmfam=\seveneufm \scriptscriptfont\eufmfam=\fiveeufm % 9

\font\teneufb=eufb10 \font\seveneufb=eufb7 \font\fiveeufb=eufb5 \newfam\eufbfam \textfont\eufbfam=\teneufb \scriptfont\eufbfam=\seveneufb \scriptscriptfont\eufbfam=\fiveeufb % 10

\font\teneurm=eurm10 \font\seveneurm=eurm7 \font\fiveeurm=eurm5 \newfam\eurmfam \textfont\eurmfam=\teneurm \scriptfont\eurmfam=\seveneurm \scriptscriptfont\eurmfam=\fiveeurm % 11

\font\teneurb=eurb10 \font\seveneurb=eurb7 \font\fiveeurb=eurb5 \newfam\eurbfam \textfont\eurbfam=\teneurb \scriptfont\eurbfam=\seveneurb \scriptscriptfont\eurbfam=\fiveeurb % 12

\font\teneusm=eusm10 \font\seveneusm=eusm7 \font\fiveeusm=eusm5 \newfam\eusmfam \textfont\eusmfam=\teneusm \scriptfont\eusmfam=\seveneusm \scriptscriptfont\eusmfam=\fiveeusm % 13

\font\teneusb=eusb10 \font\seveneusb=eusb7 \font\fiveeusb=eusb5 \newfam\eusbfam \textfont\eusbfam=\teneusb \scriptfont\eusbfam=\seveneusb \scriptscriptfont\eusbfam=\fiveeusb % 14

\font\tenss=cmss10 \font\sevenss=cmss7 \font\fivess=cmss5 \newfam\ssfam \textfont\ssfam\tenss \scriptfont\ssfam\sevenss \scriptscriptfont\ssfam\fivess % 15
\def\sf{\fam\ssfam}

\font\sevenit=cmti7 \scriptfont\itfam=\sevenit

\let\articletitle\seventeenss
\let\chaptertitle\twelvebf
\let\sectiontitle\tenbf
\let\subsectiontitle\tenbfit
\let\subsubsectiontitle\tenit
\let\contchaptertitle\tenbf % chapter titles in the table of contents
\let\contsectiontitle\tenrm % ditto for sections
\let\contsubsectiontitle\sevenrm % ditto for subsections
\let\contsubsubsectiontitle\fiverm % ditto for subsubsections
\let\parnumfont\tenrm % paragraph numbers
\let\parbackreffont\fiverm % back references at the end of paragraphs
\let\proclaimfont\tenbf
\let\prooffont\tenit
\let\mainfont\tenrm

\def\lhighlight{} % in principle, links can be highlighted, but this is usually done by the renderer
\def\suppresscomments#1#2#3{} % for sharing drafts
\def\prohibitcomments#1#2#3{\errmessage{Draft comments are not allowed in the final version}} % for final versions

\ifx\format\undefined\else\format\fi % default paper size is letter
\ifx\comment\undefined\def\comment{\prohibitcomments}\fi

\newdimen\hmargin
\newdimen\vmargin
\newdimen\plaintextwidth
\newdimen\plaintextheight
\newdimen\totalht
\newdimen\totalwd
\def\shipbox#1{%
	\totalht\ht#1
	\advance\totalht2\vmargin
	\totalwd\wd#1
	\advance\totalwd2\hmargin
	\hoffset-1in
	\advance\hoffset\hmargin
	\voffset-1in
	\advance\voffset\vmargin
	\dumpbox\totalwd\totalht{#1}}

\newtoks\cont
\def\contents{\begingroup\suppressbackreftrue\let\anchor\gobble\the\cont\endgroup}
\let\printcont\gobble
\def\contlinechapt#1#2{\printcont{#2}\smallskip\everypar{}\noindent{\contchaptertitle\llink{chapter.#1}{#2}}\par} % Link: contents to chapter
\def\contlinesect#1#2{\printcont{#1.  #2}\everypar{}\indent{\contsectiontitle\llap{#1.\enskip}\llink{section.#1}{#2}}\par} % Link: contents to section
\def\contlineunsect#1#2{\printcont{#2}\everypar{}\noindent{\contsectiontitle\llink{section.#1}{#2}}\par} % Link: contents to section
\def\contlinesubsect#1#2{\printcont{#1.  #2}\everypar{}\indent{\contsubsectiontitle{#1.\enskip}\llink{paragraph.#1}{#2}}\par} % Link: contents to subsection
\def\contlinesubsubsect#1#2{\printcont{#1.  #2}\everypar{}\indent\indent\indent{\contsubsubsectiontitle\llap{#1.\enskip}\llink{paragraph.#1}{#2}}\par} % Link: contents to subsubsection
\def\addcont#1#2#3{\append\cont{#1}\appendexpand\cont{{#2}}\append\cont{{#3}}}

\newif\ifpresec \presecfalse % true if we are in a chapter, but before its first section
\def\chapter#1\par{%
	\def\chapname{#1}%
	\parn0
	\subparn0
	\presectrue
	\numbfalse
	\addcont\contlinechapt\chapname{#1}%
	\curverb{Chapter~}%
	\assignlabel{chapter}{\chapname}%
	\tchapter#1\par}
\def\tchapter#1\par{% typeset chapter title
	\everypar{\let\beforesect\beforesection}% do not number paragraphs before the first section
	\chapbreak\bigbreak
	\centerline{\plabel\chaptertitle\anchor{chapter.#1}{#1}}% Anchor: chapter
	\nobreak\medskip
	\let\beforesect\relax % do not break pages right after a chapter title
}

\let\chapbreak\relax % \let\chapbreak\tchapbreak to start chapters on a new page

\newcount\secn \secn0
\def\section#1\par{%
	\ifx\sectionid\undefined\advance\secn1 \edef\sectionid{\the\secn}\fi
	\everypar{\numpar}% number paragraphs inside a section
	\parn0
	\subparn0
	\presecfalse
	\numbfalse
	\addcont\contlinesect\sectionid{#1}%
	\curverb{\S}%
	\assignlabel{section}{\sectionid}%
	\tsection#1\par
}
\def\unsection#1\par{%
	\everypar{\numpar}% number paragraphs inside a section
	\parn0
	\subparn0
	\presecfalse
	\numbfalse
	\addcont\contlineunsect{#1}{#1}%
	\curverb{\S}%
	\assignlabel{section}{#1}%
	\tsection#1\par
}
\def\tsection#1\par{
	\beforesect\let\beforesect\beforesection
	\typesetsection{#1}%
	\aftersection	
	\let\sectionid\undefined
}
\def\beforesection{\vskip0pt plus.3\vsize \penalty-250 \vskip0pt plus-.3\vsize \bigskip \vskip\parskip}
\def\aftersection{\nobreak\smallskip}
\def\typesetsection#1{\leftline{\sectiontitle
	\ifx\sectionid\undefined
		\plabel\anchor{section.#1}{}%
	\else
		\hbox to \parindent{\hss\plabel\anchor{section.\sectionid}{\sectionid}\enspace\hfill}%
	\fi#1}} % Anchor: section
\let\beforesect\beforesection

\def\subsection#1\par{\bigbreak\numbtrue\curverb{\S}\subsectiontitle\noindent#1\/\mainfont\par\nobreak\medskip\addcont\contlinesubsect{\the\secn.\the\parn}{#1}}
\def\subsubsection#1\par{\bigbreak\numbtrue\curverb{\S}\subsubsectiontitle\noindent#1\/\mainfont\par\nobreak\medskip\addcont\contlinesubsubsect{\the\secn.\the\parn}{#1}}

\def\sskip#1{\ifdim\lastskip<\medskipamount \removelastskip\penalty#1\medskip\fi}
\def\slug{\hbox{\kern1.5pt\vrule width2.5pt height6pt depth1.5pt\kern1.5pt}}
\newif\ifqed \newif\ifneedqed
\def\qed{\unskip\nobreak\ \slug\ifhmode\spacefactor3000 \fi\global\qedtrue}
\def\proclaim{\medbreak\atendpar{\sskip{55}}\numbtrue\gproclaim\proclaimfont} % always number theorems etc.
\def\proof{\medbreak\atendpar{\sskip{-55}}\needqedtrue\gproclaim\prooffont}
\def\xproclaim#1.{\medbreak\atendpar{\sskip{55}}{\everypar{}\noindent}{\proclaimfont#1.\enspace}\ignorespaces} % like \proclaim, no numbering
\let\abstract\xproclaim % abstracts don't get numbered

\def\ppar{\endgraf{\everypar{}\indent}} % paragraph inside a proof, not seen by \proof, has backrefs
\newtoks\atendpar % tokens to be inserted after \endgraf
\newtoks\atendbr % paragraph backreferences
\newif\ifnumb \numbfalse % number this paragraph?
\def\finishpar{\ifhmode\ifneedqed\ifqed\else\qed\fi\qedfalse\needqedfalse\fi\iflist\endlist\fi\the\atendbr\atendbr{}\endgraf\the\atendpar\atendpar{}\numbfalse\fi}
\def\endlist{\iflist\listfalse\endgraf{\parskip\smallskipamount\everypar{}\noindent}\fi} % terminate a list without starting a new paragraph
\newcount\parn
\newcount\subparn
\newif\ifparbref
\def\nextpar{\ifnumb\advance\parn1 \printlabel{advancing paragraph number to \the\parn}\def\brt{}%
	\ifparbref\edef\cseq{\csname backreference-list.paragraph.\the\secn.\the\parn\endcsname}%
	\expandafter\ifx\cseq\relax\else\edef\brt{\cseq}\printbackref{back references for paragraph.\the\secn.\the\parn: \cseq}\fi\fi
	\assignlabel{paragraph}{\the\secn.\the\parn}\fi}
\def\numpar{\ifnumb\nextpar\expand{\noexpand\typesetparnum{\the\secn.\the\parn}\noexpand\typesetpbr{\brt}}\fi}
\newdimen\pindent \pindent\parindent

\ifx\draftnum\undefined % normal numbering
\def\gproclaim#1#2.{\curverb{#2~}% #1: font, #2: prefix
	\ifnumb\nextpar\fi
	{\everypar{}\noindent}%
	\plabel#1#2%
	\ifnumb\ \anchor{paragraph.\the\secn.\the\parn}{}\the\secn.\the\parn\fi.\enspace\mainfont % Anchor: proclaim
	\ifnumb\expand{\noexpand\typesetpbr{\brt}}\fi
	\ignorespaces}
\def\typesetparnum#1{\ifnumb{\plabel\parnumfont\anchor{paragraph.#1}{}#1.\enspace}\fi} % Anchor: explicitly numbered paragraph
\def\typesetpbr#1{\ifnumb\def\brtext{#1}\ifx\brtext\empty\else\setetok\atendbr{{\parbackreffont Used in \noexpand\stripcomma\brtext.}}\fi\fi}
\else % number everything for proofreading purposes
\let\numbfalse\relax \numbtrue
\def\gproclaim#1#2.{\curverb{#2~}\medbreak\noindent#1#2.\enspace\mainfont\ignorespaces}
\inewdimen\brwidth \brwidth.6in
\def\typesetparnum#1{\ifnumb\llap{\plabel\anchor{paragraph.#1}{}\parnumfont#1\enspace}\fi} % Anchor: implicitly numbered paragraph
\def\typesetpbr#1{\ifnumb\def\brtext{#1}\ifx\brtext\empty\else
	\llap{\smash{\vtop{\everypar{}\raggedright\rightskip0pt plus 0pt \leftskip0pt plus 1fill \hsize\brwidth
	\parnumfont \strut \break % the first line contains paragraph number
	\parbackreffont\stripcomma#1}}\enspace}\fi\fi}
\fi

\def\hang{\hangindent\pindent}
\newif\iflist
\def\textindent#1{{\everypar{}\parindent\pindent\indent}\llap{#1\enspace}\listtrue\ignorespaces}

\def\li{\item{$\bullet$}}

\def\item{\endgraf\hang\textindent}
\def\filbreak{\endgraf\vfil\penalty-200\vfilneg}

\def\eject{\endgraf\break}
\def\supereject{\endgraf\penalty-20000}
\def\smallbreak{\endgraf\ifdim\lastskip<\smallskipamount
	\removelastskip\penalty-50\smallskip\fi}
\def\medbreak{\endgraf\ifdim\lastskip<\medskipamount
	\removelastskip\penalty-100\medskip\fi}
\def\bigbreak{\endgraf\ifdim\lastskip<\bigskipamount
	\removelastskip\penalty-200\bigskip\fi}

\let\printbackref\gobble
\def\predefbackref#1{%
	\printbackref{defining back reference list backreference-list.#1}%
	\expandafter\gdef\expandafter\cseq\expandafter{\csname backreference-list.#1\endcsname}
	\expandafter\ifx\cseq\relax\expandafter\gdef\cseq{}\else\printbackref{duplicate omitted}\fi
}
\newcount\backref \backref0
\newif\ifsuppressbackref \suppressbackreffalse
\def\firstletter#1#2\endletter{#1}
\newtoks\backreflist
\newif\ifaddbr \addbrfalse
\def\recordbackref#1{%
	\edef\params{{\ifpresec\expandafter\firstletter\chapname\endletter\else\the\secn\fi.\the\parn\ifnumb\else*\fi}{\the\backref}{\the\inputlineno}}%
	\printbackref{recording back reference \string#1 for future processing with params \params}%
	\edef\tmp{\the\backreflist\noexpand\processbackref\noexpand#1\params}%
	\global\backreflist\expandafter{\tmp}%
	\printbackref{new content of backreflist: \the\backreflist}} % \global is needed because this can be invoked inside {\it ...}, say
\def\processbackref#1#2#3#4{%
	\edef\key{\expandafter\gobble\string#1}%
	\edef\cseq{\csname id.\key\endcsname}%
	\expandafter\ifx\cseq\relax \ewarningline{#4}{Undefined reference \string#1.}{\string#1}\else
		\edef\lseq{\csname backreference-list.\cseq\endcsname}%
		\printbackref{processing back reference number #3 \string#1, originating from #2, at line #4; adding to \lseq}%
		\edef\lastnumber{\csname lastnumber.\cseq\endcsname}%
		\edef\newnumber{#2}%
		\addbrtrue
		\printbackref{lastnumber: \lastnumber; newnumber: \newnumber;}%
		\expandafter\ifx\csname lastnumber.\cseq\endcsname\relax % we are the first back reference
		\else\ifx\lastnumber\newnumber % same as the last one
			\printbackref{Suppressing duplicate back reference #2.}%
			\addbrfalse
		\fi\fi
		\expandafter\edef\csname lastnumber.\cseq\endcsname{#2}% record the new number
		\ifaddbr
			\expandafter\expandafter\expandafter\gdef\expandafter\expandafter\csname backreference-list.\cseq\endcsname\expandafter{\lseq, \llink{backreference.#3}{#2}}% Link: from a back reference list to the point of origin
		\fi
	\fi
}
\newtoks\labelinitlist
\def\xxstripcomma, {}
\def\xstripcomma{\if\ntok,\let\xcont\xxstripcomma\else\let\xcont\relax\fi\xcont}
\def\stripcomma{\futurelet\ntok\xstripcomma}

\inewif\ifrecorddups \recorddupstrue
\def\checkduplicates#1#2{\edef\key{\expandafter\gobble\string#1}%
	\iftypesetting\else
	\expandafter\ifx\csname line:\key\endcsname\relax\printlabel{keydefline: relax}\ifrecorddups\expandafter\xdef\csname line:\key\endcsname{\the\inputlineno}\fi\else\edef\keydefline{\csname line:\key\endcsname}\errmessage{#2}\fi\fi}

\let\printlabel\gobble
\newtoks\curverb % curverb stores the current block name, such as "Chapter", "Theorem", etc.
\ifx\draftlabel\undefined
\def\plabel{}
\else
\def\labeltext{} % label text for proofreading purposes
\def\plabel{\ifx\labeltext\empty\else\smash{\llap{\parbackreffont\labeltext\quad}}\gdef\labeltext{}\fi} % gdef because \plabel is used inside boxes
\fi
\def\assignlabel#1#2{% #1: name (e.g., reference, paragraph, section, chapter), #2: text (e.g., 2.1)
	\ifx\lastlabel\undefined\else
	\edef\key{\expandafter\expandafter\expandafter\gobble\expandafter\string\lastlabel}% label name without backslash
	\printlabel{label \key: id.\key\space = #1.#2, text.\key\space = #2}%
	\expandafter\xdef\csname id.\key\endcsname{#1.#2}% id for DVI hrefs and backreference lists
	\expandafter\xdef\csname text.\key\endcsname{#2}% text that is actually typeset
	\edef\tmp{\the\curverb}%
	\ifx\tmp\empty\else % if curverb is nonempty, define a "verbal" label with a prefix "v"
	\printlabel{label v\key: id.v\key\space = #1.#2, text.v\key\space = \the\curverb#2}%
	\expandafter\xdef\csname id.v\key\endcsname{#1.#2}% id for DVI hrefs and backreference lists
	\expandafter\xdef\csname text.v\key\endcsname{\the\curverb#2}% text that is actually typeset
	\fi
	\predefbackref{#1.#2}%
	\fi\let\lastlabel\undefined}
\newtoks\vlist % Verification list
\def\label#1{%
	\iftypesetting\else
	\checkduplicates#1{Label \string#1 was already defined at line \keydefline}%
	\addverunused#1\verifylabel
	\def\lastlabel{#1}%
	\fi
	\numbtrue % always number labeled paragraphs
}

\newif\ifyearkey % use years as bibliographic keys?
\def\y{} % use in bibliography as \y{1967}
\newdimen\bibindent % the maximum width of a bibliographic key
\newtoks\bibt % token list for all bibliographic items
\def\tbib#1{% #1 = \Paper
	\checkduplicates#1{Bibliographic reference \string#1 already defined at line \keydefline}%
	\addverunused#1\verifybib
	\edef\key{\expandafter\gobble\string#1}% reference name without backslash
	\printlabel{reference \key: id.\key\space = reference.\key, text.\key\space = \key}%
	\expandafter\edef\csname id.\key\endcsname{reference.\key}% id for DVI hrefs and backreference lists
	\expandafter\edef\csname text.\key\endcsname{\key}% text that is actually typeset
	\predefbackref{reference.\key}%
	\ifyearkey\else\setbox0=\hbox{[\key]}\ifdim\bibindent<\wd0 \bibindent=\wd0 \fi \fi % \bibindent holds the maximum length of all reference [keys]
	\appendexpand\bibt{\noexpand\typesetbib\noexpand#1\src}%
	\ifyearkey\let\next\xxftbib\else\let\next\ftbibalpha\fi\next}
\def\ftbibalpha#1\par{\append\bibt{#1\par}}
\def\xxftbib{\futurelet\next\xftbib}
\def\xftbib{\if\next[\let\next\ftbibyear\else\let\next\ftbibnoyear\fi\next}
\def\ftbibyear[#1]{\edef\yearkey{#1}\expandafter\ftbibyearbis\ignorespaces}
\def\ftbibyearbis#1\par{\append\bibt{#1\par}\ftbibend}
\def\ftbibnoyear#1\par{\append\bibt{#1\par}\edef\yearkey{\extractyear#1\par}\ftbibend}
\def\ftbibend{\expandafter\ifx\csname year.\yearkey\endcsname\relax
		\edef\yearindex{0}%
	\else
		\edef\yearindex{\csname year.\yearkey\endcsname}%
	\fi
	{\count0=\yearindex
	\advance\count0 by 1 %
	\expandafter\xdef\csname year.\yearkey\endcsname{\the\count0 }%
	\xdef\alphakey{\ifcase \count0 ?\or a\or b\or c\or d\or e\or f\or g\or h\or i\or j\or k\or l\or m\or n\or o\or p\or q\or r\or s\or t\or u\or v\or w\or x\or y\or z\else .\the\count0 \fi}}%
	\printlabel{key \key, year key \yearkey, alpha key \alphakey}%
	\expandafter\edef\csname text.\key\endcsname{\noexpand\typesetyearalpha{\key}{\yearkey}{\alphakey}}%
	\printlabel{reference \key\space adjustment: id.\key\space = reference.\key, text.\key\space = \yearkey.\alphakey}%
	\setbox0=\hbox{[\yearkey.\alphakey]}\ifdim\bibindent<\wd0 \bibindent=\wd0 \fi% \bibindent holds the maximum length of all reference [keys]
}
\def\typesetyearalpha#1#2#3{%
	\edef\yearindex{\csname year.#2\endcsname}%
	\ifnum\yearindex=1 #2\else#2.#3\fi}
\def\extractyear#1\y#2#3\par{#2}
\newif\iftype
\def\typesetbib#1#2\par{\edef\key{\expandafter\gobble\string#1}%
	\edef\bibbr{\csname backreference-list.reference.\key\endcsname}%
	\typetrue\ifsuppressunusedbib\ifx\bibbr\empty\typefalse%\warning{Suppressing unused bibliography item \string#1}
	\fi\fi
	\iftype
	\noindent\hbox to \bibindent{[\anchor{reference.\key}{\csname text.\key\endcsname}]\hfil}#2% Anchor: reference
	\ifx\bibbr\empty\else\expandafter\stripcomma\bibbr.\fi
	\hangindent\bibindent\filbreak
	\fi}

\let\printverify\gobble
\def\addverunused#1#2{\appendexpand\vlist{\noexpand#2\noexpand#1{\the\inputlineno}}}
\def\verifyref#1#2#3{\printverify{verifying for #3 \string#1 (line #2)}%
	\edef\key{\expandafter\gobble\string#1}%
	\edef\cseq{\csname id.\key\endcsname}%
	\edef\tmp{\csname backreference-list.\cseq\endcsname}%
        \ifx\tmp\empty\ewarningline{#2}{#3 \string#1.}{\string#1}\fi}
\def\verifylabel#1#2{\verifyref#1{#2}{Unused label}}
\def\verifybib#1#2{\verifyref#1{#2}{Unused reference}}

\let\printurl\gobble
\newtoks\urltext
\newtoks\urlt
\newif\ifpunct
\def\urldash{-}
\def\urltilde{{\tensy^^X}} % like \sim
\def\ndash{\def\urldash{--}}
\def\http://{\hfil\penalty900\hfilneg\urltext={http://}\urlt={http:/\negthinspace/}\punctfalse\urlgrab}
\def\https://{\hfil\penalty900\hfilneg\urltext={https://}\urlt={https:/\negthinspace/}\punctfalse\urlgrab}
\def\urlgrab{\catcode`\#=11 \catcode`\&=11 \futurelet\ntok\urldispatch}
\def\urldispatch{%
	\ifx\ntok~\let\proceed\urlcont\else
	\ifcat\noexpand\ntok\space\let\proceed\urlfinish\else
	\ifcat\noexpand\ntok\relax\let\proceed\urlfinish\else
	\let\proceed\urlcont
	\fi\fi\fi\proceed}
\def\urlcont#1{\ifpunct\appendexpand\urltext\punctc\appendexpand\urlt\punctc\punctfalse\fi
	\ifx\ntok~\appendexpand\urltext{\noexpand~}\appendexpand\urlt\urltilde
	\else\if\ntok\ampersand\appendexpand\urltext{&}\appendexpand\urlt{\&}%
	\else\if\ntok\ohash\appendexpand\urltext\ohash\appendexpand\urlt\#%
	\else\if\ntok_\appendexpand\urltext_\appendexpand\urlt\_%
	\else\if\ntok-\appendexpand\urltext-\appendexpand\urlt\urldash
	\else\if\ntok.\puncttrue\def\punctc{.}%
	\else\if\ntok,\puncttrue\def\punctc{,}%
	\else\if\ntok;\puncttrue\def\punctc{;}%
	\else\appendexpand\urltext{#1}\appendexpand\urlt{#1}%
	\fi\fi\fi\fi\fi\fi\fi\fi\urlgrab}
\def\urlfinish{\catcode`\#=6 \catcode`\&=4 \hbox{\printurl{\the\urltext}\link{\the\urltext}{\the\urlt}}\ifpunct\punctc\punctfalse\fi\def\urldash{-}}
\def\idgrab{\futurelet\ntok\iddispatch}
\def\iddispatch{\ifcat\noexpand\ntok\space\let\proceed\urlfinish
		\else\if\ntok,\let\proceed\urlfinish
		\else\let\proceed\idcont
		\fi\fi\proceed}
\def\idcont#1{\ifpunct\appendexpand\urltext.\appendexpand\urlt.\punctfalse\fi
	\if\ntok.\puncttrue\def\punctc{.}%
	\else\if\ntok_\appendexpand\urltext_\appendexpand\urlt\_%
	\else\appendexpand\urltext{#1}\appendexpand\urlt{#1}%
	\fi\fi\idgrab}

\let\printgrab\gobble
\newtoks\grabname
\newtoks\grabtoks % beginning
\newtoks\grabcseq % control sequence name
\newtoks\subsuptoks % ^ or _
\newtoks\dtoks % first letter of a subscript
\newcount\grabsize
\newif\ifgrabsubscript % include subscript?
\newif\ifgrabsupscript % include superscript?
\def\grabsequence{\bgroup % \bgroup allows us to say things like $C^@op$, with "op" being a superscript
	\grabsubscripttrue\grabsupscripttrue % do grab sub/superscripts
	\grabstring}
\def\grabalpha{\bgroup % same
	\grabsubscriptfalse\grabsupscriptfalse % do not grab sub/superscripts
	\grabstring}
\def\grabingroup{\ifinfont\errmessage{Already inside a math token}\fi\append\grabtoks{\bgroup\grablink}\infonttrue}
\def\graboutgroup{\ifinfont\append\grabtoks{\endgrablink\egroup}\infontfalse\fi}
\def\grabstring#1#2#3{% #1 = descriptive name like cat, fun, trans; #2 = font like \bf, \rm, \it; #3 = postcommand like \nolimits
	\let\specialhat^
	\catcode`\^=7
	\aftergroup#3 % #3 could be \nolimits or \limits
	\ifx\specialaddon\undefined\else\expandafter\aftergroup\specialaddon\let\specialaddon\undefined\fi
	\inewif\ifdefine \inewif\ifinfont
	\printgrab{}\printgrab{grab a string of type #1, typeset using font \string#2, with postcommand \string#3}%
	\grabname{#1}\def\grabfont{#2}\grabsize0 \grabtoks={}\grabingroup \append\grabtoks{#2}\grabcseq={}%
	\futurelet\ntok\grabdeflookahead}
\def\grabdeflookahead{\if=\noexpand\ntok % @=Set creates an anchor, whereas @Set refers to it
	\definetrue\printgrab{defining}\expandafter\grabgobblefuturelet
	\else\printgrab{referencing}\definefalse\expandafter\grablookahead\fi}
\def\grabgobblefuturelet#1{\futurelet\ntok\grabtestforsilent} % gobble = and look for another =
\newif\ifsilentgrab
\def\grabtestforsilent{\if=\noexpand\ntok \silentgrabtrue \let\ncom\grabsilenteq \else \silentgrabfalse \let\ncom\grablookahead \fi \ncom}
\def\grabsilenteq={\grabfuturelet}
\def\grabfuturelet{\futurelet\ntok\grablookahead}
\def\grablookahead{\printgrab{futurelet token meaning: \meaning\ntok}%
	\let\ncom\grabfinish
	\if\bgroup\noexpand\ntok \printgrab{left brace, terminating}%
	\else \if\egroup\noexpand\ntok \printgrab{right brace, terminating}%
	\else \if\space\noexpand\ntok \printgrab{blank space, terminating}%
	\else \let\ncom\grabexamine \fi\fi\fi \ncom}
\def\grabexamine#1{\printgrab{grabexamine argument: \string#1, meaning \meaning#1}%
	\def\ncom{\grabfinish#1}%
	\ifcat$\ifcat*\string#1\fi$% is #1 not a command sequence?
		\ifcat _\noexpand#1 \ifgrabsubscript\printgrab{subscript, continuing}%
			\graboutgroup \append\grabtoks{#1}\subsuptoks{#1}\def\ncom{\grabsubsupfuturelet}%
						\else\printgrab{subscript, terminating}\fi
		\else \ifcat ^\noexpand#1 \ifgrabsupscript\printgrab{superscript, continuing}%
			\graboutgroup \append\grabtoks{#1}\subsuptoks{#1}\def\ncom{\grabsubsupfuturelet}%
						\else\printgrab{superscript, terminating}\fi
		\else \ifx \specialhat#1 \ifgrabsupscript\printgrab{specialhat superscript, continuing}%
			\graboutgroup \append\grabtoks{#1}\subsuptoks{#1}\def\ncom{\grabsubsupfuturelet}%
						\else\printgrab{superscript, terminating}\fi
		\else \ifcat\noexpand~\noexpand#1 \printgrab{active character \string#1, examining further}%
			\ifnum1=\uccode`#1 \printgrab{UTF-8 letter, continuing}%
				\advance\grabsize1 \append\grabtoks{#1}\appendexpand\grabcseq{\string#1}\def\ncom{\grabfuturelet}%
			\else
				\ifnum\the\grabsize=0 \printgrab{Nothing grabbed so far, continuing}%
					\advance\grabsize1 \append\grabtoks{#1}\appendexpand\grabcseq{\string#1}\def\ncom{\grabfuturelet}%
				\else\printgrab{Not a UTF-8 letter and not the first character in a string, terminating}%
				\fi
			\fi
		\else \ifcat a\noexpand#1 \printgrab{letter #1, continuing}%
			\advance\grabsize1 \append\grabtoks{#1}\append\grabcseq{#1}\def\ncom{\grabfuturelet}%
		\else\printgrab{nonactive character \string#1}%
			\ifnum\the\grabsize=0 \printgrab{sole argument, adding and terminating}%
				\advance\grabsize1 \append\grabtoks{#1}\append\grabcseq{#1}\def\ncom{\grabfinish}%
			\else\printgrab{terminating}\fi
		\fi\fi\fi\fi\fi
	\else \printgrab{command sequence \string#1, terminating}\fi
	\ncom}
\def\grabsubsupfuturelet{\futurelet\ntok\grabsubsuplookahead}
\newcount\dig
\newif\ifdigit
\def\grabsubsuplookahead{\printgrab{subsup futurelet token meaning: \meaning\ntok}%
	\if\bgroup\noexpand\ntok \printgrab{left brace, continuing}\let\ncom\grabentiresubsup%
	\else \if\egroup\noexpand\ntok \printgrab{right brace, continuing}\let\ncom\grabentiresubsup%
	\else \if\space\noexpand\ntok \errmessage{Blank space after \the\subsuptoks}%
	\else \let\ncom\grabsubsupexamine \fi\fi\fi \ncom}
\def\grabentiresubsup#1{\printgrab{subsup entire group added}\grabingroup\append\grabtoks{#1}\graboutgroup\grabfuturelet}
\def\grabsubsupexamine#1{\printgrab{examining subsup argument \string#1, meaning \meaning#1}%
	\ifcat$\ifcat*\string#1\fi$% is #1 not a command sequence?
		\ifcat\noexpand~\noexpand#1 \printgrab{active character \string#1, continuing}%
			\grabingroup\appendexpand\grabtoks{\grabfont\noexpand#1}\let\ncom\grabfuturelet
		\else\ifnum"8000=\the\mathcode`#1 \printgrab{math active character \string#1, continuing}%
			\let\specialaddon\egroup
			\def\ncom{#1}%
			\grabtypeset
		\else\ifcat a\noexpand#1 \printgrab{letter #1, checking whether single or not}%
			\dtoks{#1}\def\ncom{\futurelet\ntok\grabsubsupsecondletterlookahead}%
		\else\printgrab{something else, inserting a single-character sub/superscript, continuing}
			\appendexpand\grabtoks{\bgroup\grabfont\noexpand#1\egroup}%
			\advance\grabsize2 \appendexpand\grabcseq{\the\subsuptoks\string#1}% possibly ignore digits here
			\def\ncom{\grabfuturelet}\fi\fi\fi
	\else \printgrab{command sequence \string#1, continuing}%
		\append\grabtoks{#1}\let\ncom\grabfuturelet\fi
	\ncom}
\def\grabsubsupsecondletterlookahead{\def\ncom{\appendexpand\grabtoks{\the\dtoks}\grabfuturelet}%
	\ifcat a\noexpand\ntok \printgrab{not a single letter, grabbing the entire subsupscript}%
		\advance\grabsize1 \appendexpand\grabcseq{\expandafter\string\the\subsuptoks}% append _ or ^ to the label
		\grabingroup
		\appendexpand\grabtoks{\grabfont\the\dtoks}%
		\appendexpand\grabcseq{\the\dtoks}%
		\def\ncom{\grabfuturelet}%
	\else \printgrab{single letter, continuing}\fi\ncom}
\def\grabfinish{\printgrab{grabfinish}\graboutgroup\grabtypeset\egroup}
\def\grabtypeset{\printgrab{grabtypeset grabsize=\the\grabsize, grabtoks=\the\grabtoks, grabcseq=\the\grabcseq}%
	\def\grablink##1\endgrablink{##1}%
	\ifnum\the\grabsize=0 \errmessage{No string to grab}\fi
	\ifnum\the\grabsize>1 % single-letter names are not references
		\ifdefine % defining a mathematical identifier
			\expandafter\checkduplicates\csname\the\grabname.\the\grabcseq\endcsname{Mathematical identifier \key\space already defined at line \keydefline}%
			\iftypesetting % if we are actually typesetting, create an anchor
				\ifsilentgrab
					\expandafter\gdef\csname silent:\the\grabname.\the\grabcseq\endcsname{}% record that this id is silent
				\else
					\anchor{\the\grabname.\the\grabcseq}{}% Anchor: definition of a mathematical identifier
				\fi % only create if not silent
			\else
				\expandafter\xdef\csname id.\the\grabname.\the\grabcseq\endcsname{paragraph.\the\secn.\the\parn}% paragraph id for a back reference list
				\predefbackref{paragraph.\the\secn.\the\parn}%
  			\fi
		\else % referencing a mathematical identifier
			\iftypesetting % we are actually typesetting
				\global\advance\backref1
				\expandafter\ifx\csname line:\the\grabname.\the\grabcseq\endcsname\relax % identifier is undefined
					\warning{Undefined mathematical identifier \the\grabname.\the\grabcseq}%
					\expandafter\gdef\csname\the\grabname.\the\grabcseq\endcsname{\relax}% report undefined references only once; gdef because inside bgroup..egroup
				\else % identifier has ben defined
					\expandafter\ifx\csname silent:\the\grabname.\the\grabcseq\endcsname\empty
						\edef\grablink##1\endgrablink{{##1}}%
					\else % not silent, need a hyperlink
						\edef\grablink##1\endgrablink{\noexpand\llink{\the\grabname.\the\grabcseq}{##1}% Link: from a mathematical identifier to its definition 
							\noexpand\anchor{backreference.\the\backref}{}}% Anchor: back reference for a mathematical identifier reference
					\fi
				\fi
			\else % not yet typesetting, just collecting back references
				\ifsuppressbackref\else
					\global\advance\backref1
					\printbackref{math back reference to \the\grabname.\the\grabcseq: backref.\the\backref\space at line \the\inputlineno}%
					\expandafter\recordbackref\csname\the\grabname.\the\grabcseq\endcsname
				\fi
			\fi
		\fi
	\fi
	\ifsilentgrab\else\the\expandafter\grabtoks\fi}
\inewif\ifsilent \inewif\ifanchor
\newif\ifsuppresscs
\catcode`\^=7   \def\specialhat{\ifmmode\def\next{^}\else\let\next\beginxref\fi\next} \catcode`\^=\active \let^=\specialhat
\def\silentxref#1{\futurelet\next\silentxrefswitch}
\def\silentxrefswitch{\silenttrue\xref}
\def\beginxref{\futurelet\next\beginxrefswitch}
\def\beginxrefswitch{\ifx\next\specialhat\let\next\silentxref \else\silentfalse\let\next\xref\fi \next}
\def\xref{\leavevmode\futurelet\next\xrefswitch}
\def\xrefswitch{\ifx\next!\let\next\verbalxref \else \ifx\next=\let\next\anchorxref \else \anchorfalse \let\next\normalxref \fi \fi \next}
\newtoks\vtoksl
{\count0="C2 \loop\ifnum\count0<"F5 \catcode\count0=11 \advance\count0 by 1 \repeat % make all legitimate initial UTF-8 octets letters
\gdef\plainaccents{\suppresscstrue%
	\def\`##1{##1\empty ̀}%
	\def\'##1{##1\empty ́}%
	\def\^##1{##1\empty ̂}%
	\def\"##1{##1\empty ̈}%
	\def\~##1{##1\empty ̃}%
	\def\=##1{##1\empty ̄}%
	\def\.##1{##1\empty ̇}%
	\def\u##1{##1\empty ̆}%
	\def\v##1{##1\empty ̌}%
	\def\H##1{##1\empty ̋}%
	\def\t##1{##1\empty ͡}%
}}
\def\plainaccents{\let\xcsname=\empty \let\xendcsname=\empty}
\def\verbalxref!{\begingroup\plainaccents\verbalxrefaux}
\def\verbalxrefaux#1{%
	\lowercase{\vtoksl{#1}}%
	\expandtoks\vtoksl
	\iftypesetting
		\expandafter\gdef\expandafter\cseq\expandafter{\csname verbal.\the\vtoksl\endcsname}%
		\printlabel{verbal xref \the\vtoksl}%
		\expandafter\ifx\cseq\relax\errmessage{Undefined reference to \the\vtoksl}\else\cseq\fi
	\else
		\ifsuppressbackref\else
			\global\advance\backref1 %
			\printbackref{label verbal.#1: backref.\the\backref\space at line \the\inputlineno}%
			\blah
			\expandafter\recordbackref\csname verbal.\the\vtoksl\endcsname
		\fi
	\fi
	\endgroup}
\def\initverballabelcommand#1{%
	\printlabel{initializing verbal label #1 (\the\curverb\the\secn.\the\parn)}%
	\expandafter\xdef\csname id.verbal.#1\endcsname{paragraph.\the\secn.\the\parn}%
	\ifnum\parn=0 %
	\expandafter\xdef\csname text.verbal.#1\endcsname{\the\curverb\the\secn}%
	\else
	\expandafter\xdef\csname text.verbal.#1\endcsname{\the\curverb\the\secn.\the\parn}%
	\fi
	\expandafter\initlabelcommand\csname verbal.#1\endcsname
}
\def\anchorxref={\anchortrue\futurelet\next\anchorxrefswitch}
\def\anchorxrefswitch{\ifx\next:\let\next\nonitalicanchor\else\italictrue\let\next\normalxref\fi \next}
\def\nonitalicanchor:{\italicfalse\normalxref}
\newtoks\firsttoks \newtoks\secondtoks \inewif\ifplural \inewif\ifitalic
\def\parseplural#1[#2|#3]{\let\next\parseplural\ifx\hfuzz#2\hfuzz\ifx\hfuzz#3\hfuzz\let\next\relax\else\pluraltrue\fi\else\pluraltrue\fi
	\append\firsttoks{#1#2}\append\secondtoks{#1#3}\next}
\newtoks\firsttoksl
\newtoks\secondtoksl
\newtoks\nexttoks
\def\normalxref{\begingroup\plainaccents\normalxrefaux}
\def\normalxrefaux#1{\firsttoks{}\secondtoks{}\pluralfalse\parseplural#1[|]%
	\lowercase\expandafter{\expandafter\firsttoksl\expandafter{\the\firsttoks}}%
	\expandtoks\firsttoksl
	\lowercase\expandafter{\expandafter\secondtoksl\expandafter{\the\secondtoks}}%
	\expandtoks\secondtoksl
	\ifanchor
		\iftypesetting\else
			\initverballabelcommand{\the\firsttoksl}%
			\ifplural\initverballabelcommand{\the\secondtoksl}\fi
			\predefbackref{paragraph.\the\secn.\the\parn}%
			\expandafter\checkduplicates\csname\the\firsttoksl\endcsname{Verbal label \the\firsttoksl\space was already defined at line \keydefline}
			\ifplural\expandafter\checkduplicates\csname\the\secondtoksl\endcsname{Verbal label \the\secondtoksl\space was already defined at line \keydefline}\fi
		\fi
		\anchor{verbal.\the\firsttoksl}{}% Anchor: verbal reference
		\ifplural\anchor{verbal.\the\secondtoksl}{}\fi % Anchor: plural verbal reference
	\else
		\ifsuppressbackref\else
			\global\advance\backref1
			\printbackref{label normal.#1: backref.\the\backref\space at line \the\inputlineno}%
			\iftypesetting
			\else
				\expandafter\recordbackref\csname verbal.\the\secondtoksl\endcsname
			\fi
			\anchor{backreference.\the\backref}{}% Anchor: back reference for a verbal reference
		\fi
	\fi
	\nexttoks{}%
	\ifsilent
		\nexttoks{\ignorespaces}%
	\else
		\ifanchor
			\ifitalic\nexttoks\expandafter{\expandafter\bgroup\expandafter\it\the\firsttoks\italcorr}%
			\else\nexttoks\expandafter{\the\firsttoks}%
			\fi
		\else
			\edef\tmp{{verbal.\the\secondtoksl}}%
			\nexttoks\expandafter{\expandafter\llink\tmp}% Link: from a verbal reference to its definition
			\nexttoks\expandafter\expandafter\expandafter{\expandafter\the\expandafter\nexttoks\expandafter{\the\firsttoks}}%
		\fi
	\fi
	\expandafter\endgroup\the\nexttoks}
\def\italcorr{\futurelet\next\italcorrtest}
\def\italcorrtest{\if,\noexpand\next\else\if.\noexpand\next\else\/\fi\fi\egroup}

\def\numdam{\urltext={http://www.numdam.org/item/?id=}\urlt={numdam:}\punctfalse\idgrab}

\def\gen:{\http://libgen.rs/book/index.php?md5=}
\def\jstor:{\https://www.jstor.org/stable/}
\def\eudml:{\https://eudml.org/doc/}
\def\repo#1#2{\urltext={#1}\urlt={#2}\punctfalse\idgrab}
\def\arXiv:{\urltext={https://arxiv.org/abs/}\urlt={arXiv:}\punctfalse\idgrab}

\def\Zbl:{\urltext={https://zbmath.org/?q=an:}\urlt={Zbl:}\punctfalse\idgrab}
\def\doi:{\ndash\urltext={https://doi.org/}\urlt={doi:}\punctfalse\urlgrab}

\def\sqr#1#2{{\thinspace\vbox{\hrule height.#2pt \hbox{\vrule width.#2pt height#1pt \kern#1pt \vrule width.#2pt} \hrule height0pt depth.#2pt}\thinspace}}
\def\square{\mathchoice\sqr64\sqr64\sqr{4.2}3\sqr33}

\def\ltoarr#1{\mathop{\count0=#1 \loop\ifnum\count0>0 \smash-\mkern-7mu \advance\count0 -1 \repeat \mathord\rightarrow}\limits} % parametrized \rightarrowfill
\def\lto#1#2{\mathrel{\ltoarr{#1}^{#2}}} % parametrized \rightarrowfill, with a label
\def\longto#1^#2_#3{\mathrel{\ltoarr{#1}^{#2}_{#3}}} % parametrized \rightarrowfill, with a label above and below
\def\lgetsarr#1{\mathop{\mathord\leftarrow \count0=#1 \loop\ifnum\count0>0 \mkern-7mu\smash-\advance\count0 -1 \repeat}\limits} % parametrized \leftarrowfill
\def\longgets#1^#2_#3{\mathrel{\lgetsarr{#1}\limits^{#2}_{#3}}} % parametrized \leftarrowfill, with a label

\def\toto{\mathrel{\vcenter{\hbox{$\to$}\kern-1.5ex \hbox{$\to$}}}}
\def\prearrfill{\smash-\mkern-7mu}
\def\postarrfill{\mkern-7mu\smash-}
\def\midarrfill#1{\cleaders\hbox{$\mkern-2mu\smash-\mkern-2mu$}\hskip0pt plus #1fil}
\def\rightarrfill{\mkern-7mu\mathord\rightarrow}
\def\leftarrfill{\mathord\leftarrow\mkern-7mu\midarrfill1\postarrfill}
\def\ltoto#1#2#3{\ifinner
	\mathrel{\vcenter{\hbox to #1em{$\prearrfill\midarrfill1{\scriptstyle#2}\midarrfill3 \rightarrfill$}%
		\kern-1.5ex \hbox to #1em{$\prearrfill\midarrfill3{\scriptstyle#3}\midarrfill1 \rightarrfill$}}}%
	\else
	\mathrel{\mathop{\vcenter{\hbox to #1em{\rightarrowfill}%
		\kern-1.5ex \hbox to #1em{\rightarrowfill}}}\limits^{#2}_{#3}}%
	\fi}
\def\ltogets#1#2#3{\ifinner
	\mathrel{\vcenter{\hbox to #1em{$\prearrfill\midarrfill1{\scriptstyle#2}\midarrfill3 \rightarrfill$}%
		\kern-1.5ex \hbox to #1em{$\leftarrfill\midarrfill3{\scriptstyle#3}\midarrfill1 \postarrfill$}}}%
	\else
	\mathrel{\mathop{\vcenter{\hbox to #1em{\rightarrowfill}%
		\kern-1.5ex \hbox to #1em{\leftarrowfill}}}\limits^{#2}_{#3}}%
	\fi}

\def\ltogetscore#1#2{\dimen0=\fontdimen6 #1 2 \divide\dimen0 by 2 \multiply\dimen0 by #2 \vcenter{\hbox to \dimen0{\rightarrowfill}\kern-1.8ex \hbox to \dimen0{\leftarrowfill}}}
\def\ltogets#1#2#3{\mathrel{\mathop{\mathchoice{\ltogetscore\textfont{#1}}{\ltogetscore\textfont{#1}}{\ltogetscore\scriptfont{#1}}{\ltogetscore\scriptscriptfont{#1}}}^{#2}_{#3}}}

\def\rx#1#2{\rlap{\kern #1pt \raise#1pt \hbox{#2}}}
\def\dottednearrow{\rx{-8}. \rx{-6}. \rx{-4}. \rx{-2}. \rx0. \rx2. \rx4. \kern6pt \raise7.7pt \hbox{$\nearrow$}}

\def\gmatrix#1#2{\null\,\vcenter{\normalbaselines
	\ialign{#1\crcr
		\mathstrut\crcr\noalign{\kern-\baselineskip}
		#2\crcr\mathstrut\crcr\noalign{\kern-\baselineskip}}}\,}
\def\cdmatrix{\gmatrix{\hfil$##$\hfil&&\enspace\hfil$##$\hfil\enspace&\hfil$##$\hfil}}
\def\sqmatrix{\gmatrix{\hfil$##$&\enspace\hfil$##$\hfil\enspace&$##$\hfil}}
\def\cdbl{\def\normalbaselines{\baselineskip20pt \lineskip3pt \lineskiplimit3pt }}

\def\cd{\cdbl\cdmatrix}
\def\sqcd{\cdbl\let\vagap\;\sqmatrix}

\newcount\arrowsize \arrowsize3
\def\mapright#1{\smash{\lto\arrowsize{#1}}}

\def\rvagap{\vagap} \def\lvagap{\vagap} \def\rvaskip{\vaskip} \def\lvaskip{\vaskip} \def\vaskip{} \def\vagap{}
\def\mapdown#1{\rvagap\Big\downarrow\rlap{$\vcenter{\hbox{$\scriptstyle#1$}}$}\rvaskip}

\def\lmapdown#1{\lvaskip\llap{$\vcenter{\hbox{$\scriptstyle#1$}}$}\Big\downarrow\lvagap}

\newcount\forno \forno0
\def\arrno#1#2{\global\advance\arr1 \edef\eeqnno{\the\arr}%
	\global\advance\forno1 \edef\eforno{\the\forno}%
	\xdef#2{\noexpand\llink{equation.\eforno}{\eeqnno}}% Link: from an equation number to its definition; xdef to ensure it can be seen outside of \wrapdiagram
	#1{(\anchor{equation.\eforno}{\eeqnno})}} % Anchor: equation number
\newbox\mdiag
\def\wrapdiagram{%
	\setbox\mdiag\vtop\bgroup
	\null % an empty hbox to ensure proper vertical positioning
	\vskip\baselineskip
	\inewcount\arr \arr0
	\baselineskip0pt
	\lineskip4pt
	\lineskiplimit4pt
	\let\par\cr
	\obeylines
	\halign\bgroup\hfil$\displaystyle##$\hfil\cr
	\ewrapdiagram}

\def\ewrapdiagram#1{#1
	\egroup
	\egroup
	\vskip0pt plus \dp\mdiag \penalty-250 \vskip0pt plus-\dp\mdiag % ensure the diagram fits on a single page
	\hangafter-\dp\mdiag
	\divide\hangafter\baselineskip
	\advance\hangafter-2
	\hangindent-\wd\mdiag
	\advance\hangindent-2em
	\hbox to\hsize{\hfil\dp\mdiag0pt \box\mdiag}%
	\ignorespaces}

\inewtoks\preamble
{\endlinechar=10 \catcode`#=12 \global\preamble{
prologues := 3;

verbatimtex
\let\endprolog

\expandafter\gobbleinit\input }\global\appendonceexpand\preamble\jobname\global\append\preamble{
\catcode`\^=7
etex

input cmarrows
setup_cmarrows(arrow_name = "texarrow"; parameter_file = "cmr10.mf"; macro_name = "drawarrow");
setup_cmarrows(arrow_name = "doublearrow"; parameter_file = "cmr10.mf"; macro_name = "drawdarrow");
def drawmarrow expr p = _apth:=p; _finmarr enddef;
rule_thickness#:=.4pt#;    % cmr10.mf: thickness of lines in math symbols
def _finmarr text t_ =
  drawarrow subpath(0, 0.5 * length(_apth)) of _apth t_;
  draw subpath(0.5 * length(_apth), length(_apth)) of _apth withpen pencircle scaled rule_thickness# t_;
enddef;

def object(suffix O)(expr x,y)(expr l) =
  save O;       
  pair O;
  O := (x,y) * u;
  picture O.tx;    
  O.tx := thelabel(l,O);
  draw O.tx;                         
enddef;                                    

def smorphism(suffix A,B) =
  save ss, tt;
  ss := xpart ((A..B) intersectiontimes bbox A.tx);
  tt := xpart ((A..B) intersectiontimes bbox B.tx);
  drawarrow subpath(ss,tt) of (A..B);
enddef;

def morphism(suffix A,B)(expr l)(expr f)(suffix $) =
  smorphism(A, B);
  label.$(l, point f[ss, tt] of (A..B));
enddef;
}}%

\newcount\vertex                                      
\ifx\mathspecials\undefined\def\mathspecials{}\fi
\def\diagram{\begininlinemp\buffertoks={}\grabdiagram}
\def\grabdiagram{\mathspecials\grabdiagramaux}                                                 
\def\grabdiagramaux#1;{\vertex0 \dcom#1,*}                          
\def\dcom{\futurelet\next\dcomswitch}                                          
\def\dcomswitch{\ifcat a\noexpand\next \let\next\dcomalpha                                       
        \else\if*\noexpand\next \addcode{\the\buffertoks}\endinlinemp\let\next\gobble                                                 
        \else\if[\noexpand\next \let\next\grabscale                                                                
        \else\errmessage{Unrecognized diagram command \next}\let\next\relax\fi\fi\fi\next}
\def\grabscale[#1]{\dimen0=#1 \addcode{save u; u = \the\dimen0;}\dcom}
\def\dcomalpha#1{\def\objectname{#1}\futurelet\next\dcomalphaswitch}
\def\dcomalphaswitch{\ifcat a\noexpand\next \let\next\dcommorphism\else                                                          
        \expandafter\edef\csname vertex:\objectname\endcsname{\the\vertex}%              
        \if:\noexpand\next\let\next\grabcoords                                              
        \else\if=\noexpand\next\let\next\grabobjectequ                                     
        \else\errmessage{Expected : or = while processing a diagram object, got \meaning\next}\let\next\relax\fi\fi\fi\next}
\def\grabcoords:#1,#2={\toks\vertex{#1,#2}\grabobjectlabel}
\def\grabobjectequ={\edef\tmp{\the\toks\vertex}\ifx\tmp\empty\errmessage{No coordinates specified for vertex \the\vertex: \objectname}\fi\grabobjectlabel}
\def\grabobjectlabel#1,{\addcode{object(\objectname, \the\toks\vertex, }\addunexpandedcode{btex #1 etex);}\advance\vertex1 \dcom}
\def\dcommorphism#1{\def\tobjectname{#1}%
        \def\labelpos{.5}%                                      
        \futurelet\next\dcommorphismswitch}                     
\def\dcommorphismswitch{\if.\noexpand\next \expandafter\grabmorphismdir \else \setupmorphismlabel \expandafter\grabmorphismpos\fi}
\def\setupmorphismlabel{
        \edef\vlabeldira{direction.\objectname.\tobjectname}%
        \edef\vlabeldirb{direction.\tobjectname.\objectname}%
        \expandafter\ifx\csname\vlabeldira\endcsname\relax
        	\expandafter\ifx\csname\vlabeldirb\endcsname\relax
		        \edef\tmpa{\csname vertex:\objectname\endcsname}\expandafter\ifx\tmpa\relax\errmessage{No such vertex: \objectname}\fi
		        \edef\tmpb{\csname vertex:\tobjectname\endcsname}\expandafter\ifx\tmpb\relax\errmessage{No such vertex: \tobjectname}\fi
		        \ifnum \tmpa<\tmpb \edef\vlabeldir{direction.\tmpa.\tmpb}\else \edef\vlabeldir{direction.\tmpb.\tmpa}\fi
		        \expandafter\ifx\csname\vlabeldir\endcsname\relax\errmessage{No label direction specified for morphism \objectname->\tobjectname}\fi
		        \edef\labeldir{\csname\vlabeldir\endcsname}%                        
		\else
        		\edef\labeldir{\csname\vlabeldirb\endcsname}%                        
		\fi
	\else
        	\edef\labeldir{\csname\vlabeldira\endcsname}%                        
	\fi
}                                                                                        
\def\grabmorphismdir.{\def\labeldir{}\futurelet\next\grabmorphismdirec}
{\catcode`\@=13                                                               
\gdef\grabmorphismdirec{\ifx@\next \let\next\grabmorphismposition
        \else\if=\noexpand\next \let\next\grabmorphismequ               
        \else \let\next\grabmorphismdirect\fi\fi\next}                       
\gdef\grabmorphismdirect#1{\edef\labeldir{\labeldir#1}\futurelet\next\grabmorphismdirec}
\gdef\grabmorphismpos{\ifx@\next \expandafter\grabmorphismposition \else \expandafter\grabmorphismequ\fi} 
\gdef\grabmorphismposition@#1={\def\labelpos{#1}\grabmorphismlabel}
}
\def\grabmorphismequ={\grabmorphismlabel}
\def\grabmorphismlabel#1,{%
        \addcodebuffer{morphism(\objectname, \tobjectname, }%
        \addunexpandedcodebuffer{btex \everymath{\scriptstyle}#1 etex, }%
        \addcodebuffer{\labelpos, \labeldir);}%
        \dcom}

\let\printextract\gobble
\def\hinitlabelcommand#1{\printextract{initializing label command \string#1}%
	\gdef#1{\printlabel{invoked label \string#1}% #1: label command, #2: label id for DVI, #3: typesetted text
		\ifsuppressbackref
			\edef\key{\expandafter\gobble\string#1}%
			\llink{\csname id.\key\endcsname}{\csname text.\key\endcsname}% Link: from a label (theorem or bibliography) to its definition
		\else
			\global\advance\backref1 %
			\printbackref{label \expandafter\gobble\string#1: backref.\the\backref\space at line \the\inputlineno}%
			\iftypesetting
			\else
				\blah
				\recordbackref#1
			\fi
			\anchor{backreference.\the\backref}{}% Anchor: back reference for a label (theorem or bibliography)
			\edef\key{\expandafter\gobble\string#1}%
			\llink{\csname id.\key\endcsname}{\csname text.\key\endcsname}% Link: from a label (theorem or bibliography) to its definition
		\fi
		}}
\let\initlabelcommand\hinitlabelcommand
\def\pinitlabelcommand#1{\printextract{initializing label command \string#1}%
	\gdef#1{\printlabel{invoked plain label \string#1}% #1: label command, #2: label id for DVI, #3: typesetted text
		\iftypesetting
			\edef\key{\expandafter\gobble\string#1}%
			\csname text.\key\endcsname
		\else
			\blah
		\fi}}
\newread\labelin
\newif\iflabelcont
\let\terminate=\relax % allow \terminate to be read by \read
\long\def\labelauxaux#1\terminate
{}
\def\preprocesslabel#1{%
	\printextract{Processing label \string#1}%
	\edef\key{\expandafter\gobble\string#1}%
	\initlabelcommand#1
	\expandafter\initlabelcommand\csname v\key\endcsname
	\labelauxaux
}
\def\preprocessbib#1{%
	\printextract{Processing bib \string#1}%
	\initlabelcommand#1
	\labelauxaux
}
\long\def\labelaux#1{\ifx#1\label\let\next\preprocesslabel\else\ifx#1\bib\let\next\preprocessbib\else\let\next\labelauxaux\fi\fi\next}
\def\processoneline{\expandafter\labelaux\labelline\relax\relax\relax\terminate}
\def\preprocess#1{% preprocessing stage to define all labels, back references, bibliographic items, and the table of contents
	\openin\labelin=#1 \labelconttrue
	\loop
		\read\labelin to\labelline
		\ifeof\labelin\let\next\labelcontfalse\else\let\next\processoneline\fi
		\next
	\iflabelcont\repeat
	\expandafter\gobbleinit\input#1
	\ifhmode\par\fi\vfill\eject % for LaTeX
	\ifhmode\par\fi\vfill\supereject
}

\def\importlabels#1{ % import labels from an external document
	\recorddupsfalse
	\let\initlabelcommand\pinitlabelcommand
	\preprocess{#1}%
	\let\initlabelcommand\hinitlabelcommand
	\recorddupstrue
	\cont={}% discard the external table of contents
	\vlist={}% discard the verification list
	\let\chapname\undefined \secn0 \parn0 \backref0 % reset the numbering
}

\newif\iftypesetting % are we in the final typesetting stage?

\typesettingfalse % no typesetting at this stage
\def\blah{blah} % placeholder for future references
\output{\setbox0\box255 \setbox0\box\footins \deadcycles0 } % fake typesetting
\let\bib\tbib % collect bibliographic items at this stage
\let\refs\relax % do not typeset bibliography
\everypar{\numpar}
\def\par{\finishpar} % no paragraph back references
\newif\ifmetapost \metapostfalse
\edef\filestem{\jobname.gen}
\newcount\figno
\newwrite\mpout
\newtoks\outtoks
\def\inimp{\ifmetapost\else\global\metaposttrue\immediate\openout\mpout=\filestem.mp \addcode{\the\preamble}\fi}
\newdimen\oldhsize
\oldhsize=\hsize
\newdimen\oldvsize
\oldvsize=\vsize
\hsize\maxdimen
\vsize\maxdimen
\preprocess\jobname % collect all labels and bibliography from the main documents and possibly also external labels
\hsize\oldhsize
\vsize\oldvsize
\ifx\fmtname\plainfmtname
\hsize\plaintextwidth
\vsize\plaintextheight
\fi
\printerr{(\jobname.tex}
\printbackref{list of all back references: \the\backreflist}%
\the\backreflist % process all back references
\ifsuppressunusedbib\def\verifybib#1#2{}\fi
\the\vlist % verify if there are any unused references or bibliography items
\printerr{)}
\ifmetapost
\addcode{end}
\immediate\closeout\mpout
\immediate\write18{mpost -interaction nonstopmode \filestem.mp}% can also compile the METAPOST file separately
\ifeof18 \warning{Compile the METAPOST file \filestem.mp manually using mpost \filestem.mp}\let\runmp\warning\fi
\fi
\let\addcode\gobble % do not write to the closed file again

\typesettingtrue
\def\importlabels#1{} % stop collecting labels and bibliographic items from other documents
\ifscroll \vsize\maxdimen\inewbox\abox\output{\setbox\abox\vbox{\unvbox255\unskip}\shipbox\abox}%
\else % real typesetting
	\ifx\fmtname\plainfmtname
	\totalht\plaintextheight
	\advance\totalht.1in
	\advance\totalht2\vmargin
	\totalwd\plaintextwidth
	\advance\totalwd2\hmargin
	\setpapersize\totalwd\totalht
	\advance\hoffset\hmargin
	\advance\voffset\vmargin
	\fi
	\output{\plainoutput}
\fi
\let\blah\undefined % disallow placeholder references
\def\bib#1\par{} % do not collect bibliographic items at this stage
\def\refs{\raggedright\rightskip0em plus \maxdimen \advance\bibindent1em \everypar{}\the\bibt\ignorespaces} % typeset bibliography at this stage
\def\addverunused#1#2{} % do not collect verification data at this stage
\def\addcont#1#2#3{} % do not touch the table of contents at this stage
\def\prepass#1{} % do not collect labels and bibliography from other files at this stage
\parbreftrue % typeset paragraph back references at this stage
\everypar{\numpar}
\let\chapname\undefined \secn0 \backref0 % reset the numbering
\expandafter\gobbleinit\input\jobname\relax
\ifhmode\par\fi\vfill\eject % for LaTeX
\ifhmode\par\fi\vfill\supereject
\end

\endprolog

\mathchardef\rightleftarrows="181D
\mathchardef\boxdot"1800
\mathcode`\:="603A
\expandafter\def\csname\string∞\endcsname{\ifmmode\infty\else{\tensy1}\fi}
\expandafter\def\csname\string⇄\endcsname{\rightleftarrows}
\expandafter\def\csname\string⊡\endcsname{\boxdot}
\chardef\eogonek"A6
\font\tenecrm=ecrm10
\expandafter\def\csname\stringę\endcsname{{\tenecrm\eogonek}}

\inewtoks\title
\inewtoks\author
\title{The enriched Thomason model structure on 2-categories}
\author{Dmitri Pavlov}
\metadata{\the\title}{\the\author}

{\articletitle \leftskip0pt plus 20em \rightskip0pt plus 20em \parfillskip0pt \parindent0pt \baselineskip20pt \the\title \bigskip}
{\tabskip0pt plus 1fil \let\par\cr \obeylines \halign to\hsize{\hfil#\hfil
\chaptertitle \the\author\vadjust{\medskip}
Department of Mathematics and Statistics, Texas Tech University
\https://dmitripavlov.org/
}}

\ifx\documentclass\undefined\else
\documentclass{elsarticle}
\usepackage[T1]{fontenc}
\usepackage[margin=1in,letterpaper]{geometry}
\usepackage{silence}
\usepackage{newunicodechar}

\newunicodechar{∞}{\ifmmode\infty\else$\infty$\fi}
\newunicodechar{⇄}{\rightleftarrows}
\newunicodechar{→}{\to}
\newunicodechar{←}{\gets}
\newunicodechar{⊡}{\boxdot}
\newunicodechar{⊠}{\boxtimes}
\newunicodechar{‖}{\|}
\newunicodechar{…}{{\ifmmode\ldots\else\dots\fi}}

\mathcode`\:="603A % colon is a punctuation mark, not a binary operation
\mathchardef\colon="303A % \colon is a binary operation, not a punctuation mark

\DeclareMathAlphabet{\matheurm}{U}{eur}{m}{n}
\ifpdf
\usepackage[hidelinks]{hyperref}
\input supp-pdf.mkii

\else
\usepackage[hidelinks,hypertex]{hyperref}
\input epsf

\fi
\usepackage{amsthm,amsmath,amssymb}
\usepackage[capitalize]{cleveref}

\crefformat{equation}{(#2#1#3)}
\crefname{figure}{Figure}{Figures}
\crefformat{section}{\S#2#1#3}
\crefmultiformat{section}{\S\S#2#1#3}{ and~#2#1#3}{, #2#1#3}{, and~#2#1#3}

\def\hhyperref{\ifmmode\undefined\fi\hyperref}

\counterwithin{equation}{section}

\theoremstyle{definition}

\numberwithin{equation}{section}

\def\li{\plainitem{$\bullet$}}
\def\plainitem{\endgraf\plainhang\textindent}
\def\plainhang{\hangindent\parindent}
\def\textindent#1{\indent\llap{#1\enspace}\ignorespaces}

\def\plainmath{\mathsurround0pt }
\def\eqalign#1{\null\,\vcenter{\openup\jot\plainmath
  \ialign{\strut\hfil$\displaystyle{##}$&$\displaystyle{{}##}$\hfil
      \crcr#1\crcr}}\,}

\def\ltoarr#1{\mathop{\count0=#1 \loop\ifnum\count0>0 \smash-\mkern-7mu \advance\count0 -1 \repeat \mathord\rightarrow}\limits} % parametrized \rightarrowfill
\def\lto#1#2{\mathrel{\ltoarr{#1}^{#2}}} % parametrized \rightarrowfill, with a label
\def\longto#1^#2_#3{\mathrel{\ltoarr{#1}^{#2}_{#3}}} % parametrized \rightarrowfill, with a label above and below
\def\lgetsarr#1{\mathop{\mathord\leftarrow \count0=#1 \loop\ifnum\count0>0 \mkern-7mu\smash-\advance\count0 -1 \repeat}\limits} % parametrized \leftarrowfill
\def\longgets#1^#2_#3{\mathrel{\lgetsarr{#1}\limits^{#2}_{#3}}} % parametrized \leftarrowfill, with a label

\def\gmatrix#1#2{\null\,\vcenter{\normalbaselines
        \ialign{#1\crcr
                \mathstrut\crcr\noalign{\kern-\baselineskip}
                #2\crcr\mathstrut\crcr\noalign{\kern-\baselineskip}}}\,}
\def\cdmatrix{\gmatrix{\hfil$##$\hfil&&\enspace\hfil$##$\hfil\enspace&\hfil$##$\hfil}}
\def\sqmatrix{\gmatrix{\hfil$##$&\enspace\hfil$##$\hfil\enspace&$##$\hfil}}
\def\cdbl{\def\normalbaselines{\baselineskip20pt \lineskip3pt \lineskiplimit3pt }}

\def\cd{\cdbl\cdmatrix}
\def\sqcd{\cdbl\let\vagap\;\sqmatrix}

\inewcount\arrowsize \arrowsize3
\def\mapright#1{\smash{\lto\arrowsize{#1}}}

\def\rvagap{\vagap} \def\lvagap{\vagap} \def\rvaskip{\vaskip} \def\lvaskip{\vaskip} \def\vaskip{} \def\vagap{}
\def\mapdown#1{\rvagap\Big\downarrow\rlap{$\vcenter{\hbox{$\scriptstyle#1$}}$}\rvaskip}

\def\lmapdown#1{\lvaskip\llap{$\vcenter{\hbox{$\scriptstyle#1$}}$}\Big\downarrow\lvagap}

\def\repo#1#2#3{\href{#1#3}{#2#3}}
\def\numdam{\repo{http://www.numdam.org/item/?id=}{numdam:}}
\def\eudml{\repo{https://eudml.org/doc/}{eudml:}}
\def\arXiv{\repo{https://arxiv.org/abs/}{arXiv:}}
\def\doi{\repo{https://doi.org/}{doi:}}
\begin
{document}
{\tabskip0pt plus 1fil \halign to\hsize{\hfil#\hfil\cr
\font\articletitle=cmss17 \articletitle The enriched Thomason model structure on 2-categories \vadjust{\bigskip}\cr
\font\chaptertitle=cmbx12 \chaptertitle Dmitri Pavlov \vadjust{\medskip}\cr
Department of Mathematics and Statistics, Texas Tech University\cr
\url{https://dmitripavlov.org/}\cr
}}
\fi

\def\id{{\rm id}}
\def\Sd{\mathop{\rm Sd}\nolimits}
\def\Ex{\mathop{\rm Ex}\nolimits}
\def\fc{{\rm c}}
\def\N{{\rm N}}
\def\Free{{\rm F}}
\def\Under{{\rm U}}
\def\FreeM{{\sf L}}
\def\UnderM{{\sf R}}

\def\Env{{\sf Env}}
\def\sSet{{\sf sSet}}
\def\Mon{{\sf Mon}}
\def\Cat{{\sf Cat}}
\def\ECat{{\sf Cat}}

\def\dom{\mathop{\sf dom}\nolimits}
\def\codom{\mathop{\sf codom}\nolimits}

\hyphenation{
Grothen-dieck
co-do-main
com-po-nent-wise
lem-ma
pre-stack
pre-stacks
pre-sheaf
pre-sheaves
man-u-script
pseudo-functor
trac-ta-ble
mon-oidal
cat-e-gor-i-cal
}

%$$

\abstract Abstract.
We prove that categories enriched in the Thomason model structure
admit a model structure that is Quillen equivalent to the Bergner model structure on simplicial categories, providing a new model for $(∞,1)$-categories.
Along the way, we construct model structures on modules and monoids in the Thomason model structure
and prove that any model structure on the category of small categories that has the same weak equivalences as the Thomason model structure is not a cartesian model structure.

\tsection Contents

\contents

\section Introduction

The Thomason model structure on the category $\Cat$ of small categories was established by Thomason [\CatCMC] and Cisinski [\Dwyer].
It is Quillen equivalent to the Kan–Quillen model structure on simplicial sets,
which means that small categories can be used as models for spaces in homotopy theory.
These models had already been proved to be very productive by the time Thomason's paper [\CatCMC] came out,
e.g., in Quillen's Theorems A and B in higher algebraic K-theory (Quillen [\HAKT])
or in computing homotopy colimits of diagrams of spaces (Thomason [\Hocolim]).

It is natural to ask whether the above Quillen equivalence has an analogue for $(∞,1)$-categories.
That is to say, does the category of small categories enriched in the Thomason model structure (i.e., strict 2-categories)
admit a model structure that is Quillen equivalent to the Bergner model structure on simplicial categories?
This paper answers this question in the affirmative.

\proclaim Theorem.
(^!{enriched Thomason model structure} and ^!{enriched Thomason Quillen equivalence}.)
The category $\ECat_\Cat$ of small ^{2-categories} admits a left proper combinatorial model structure
(the ^{enriched Thomason model structure})
whose weak equivalences are ^{Dwyer–Kan equivalences}
and acyclic fibrations are ^{Dwyer–Kan acyclic fibrations}.
The Quillen equivalence
$\fc\Sd^2⊣\Ex^2\N$
(^!{properties of the Thomason adjunction})
induces a Quillen equivalence
$$L⊣R,\qquad L:\ECat_\sSet→\ECat_\Cat,\qquad R:\ECat_\Cat→\ECat_\sSet,$$
where the right adjoint functor~$R$ applies the functor $\Ex^2\N$ to every hom-object.
That is to say, the model categories of small simplicial categories (Bergner [\SCat])
and small ^{2-categories} with the ^{enriched Thomason model structure} (^!{enriched Thomason model structure})
are Quillen equivalent.

The Quillen equivalence $L⊣R$ (^!{enriched Thomason Quillen equivalence})
immediately connects (via zigzags of Quillen equivalences)
the ^{enriched Thomason model structure} on small 2-categories
to all the other models for small $(∞,1)$-categories
such as quasicategories (Joyal [\QCat, Theorem~6.12]),
relative categories (Barwick–Kan [\RelCat, Theorem~6.1]),
Segal categories (Hirschowitz–Simpson [\Desc, Théorème~2.3]; see also Bergner [\Three, Theorems 5.1 and 7.1]),
complete Segal spaces (Rezk [\CSS, Theorem~7.2]),
simplicial categories (Bergner [\SCat, Theorem~1.1]),
marked simplicial sets (Lurie [\HTT, Proposition~3.1.3.7]), etc.,
with Quillen equivalences between them established by
Bergner [\Three, Theorem~6.3, 7.5, 8.6], Joyal–Tierney [\QCSS], Joyal [\QCSC], Barwick–Kan [\RelCat, Theorem~6.1, 6.2], Lurie [\HTT, Propositions 3.1.5.3 and 3.1.5.6].
See Bergner [\Survey, \Book] for a review of these Quillen equivalences.
The adjective “enriched Thomason” refers to $(∞,1)$-categories and distinguishes this model structure
from the Thomason model structure on 2-categories (^!{previous work}), which models $(∞,0)$-categories.
It would be interesting to see how the existing results for other definitions of $(∞,1)$-categories
(see, for example, Lurie [\HTT], Cisinski [\HCHA], Riehl–Verity [\EICT])
can be (re)formulated in the setting of 2-categories.

In support of the main theorem cited above,
we prove that the Thomason model structure enjoys a collection of properties that makes it similar to the Kan–Quillen model structure on simplicial sets,
with a prominent exception of the ^{pushout product axiom}.
The latter property immediately precludes the possibility of using existing model-categorical tools
to construct the model structure of ^!{enriched Thomason model structure},
such as the theorems of Lurie [\HTT, Proposition~A.3.2.4] and Muro [\HTEC, Theorem~1.1],
which construct a model structure on the category of small categories
enriched in a combinatorial monoidal model category
satisfying some additional conditions
(Lurie: every object is cofibrant and weak equivalences are closed under filtered colimits;
Muro: the monoid axiom holds).

\proclaim Theorem.
(^!{properties of the Thomason model structure}, ^!{properties of the Thomason adjunction}, ^!{nonexistence of cartesian structures}.)
The ^{Thomason model structure}
(^!{Thomason model structure})
on the category $\Cat$ of small categories
is a proper combinatorial
model category
that is ^{tractable} (^!{tractable}), ^{pretty small} (^!{pretty small}), ^{h-monoidal} (^!{h-monoidal}), ^{flat} (^!{flat}),
and satisfies the properties of the ^{monoid axiom} (^!{monoid axiom}) other than the nonacyclic part of the ^{pushout product axiom} (^!{pushout product axiom}).
Every model structure on the same category with the same weak equivalences is not a cartesian model structure.
Furthermore, the Quillen equivalence
$$\fc\Sd^2⊣\Ex^2\N,\qquad \fc\Sd^2:\sSet→\Cat,\qquad \Ex^2\N:\Cat→\sSet$$
satisfies the conditions of a weak monoidal Quillen equivalence in the sense of Schwede–Shipley [\EMMC, Definition~3.6],
except for the nonacyclic part of the ^{pushout product axiom}.

\subsection Previous work
^^={previous work}

Lack [\Bicat, Theorem~4]
constructs a model structure on the category of small 2-categories
whose weak equivalences are equivalences of 2-categories.
Worytkiewicz–Hess–Parent–Tonks [\TwoCatModel, \TwoCatModelCorr],
Ara–Maltsiniotis [\Vers],
Chiche [\Lax, Théorème~7.9],
Ara [\AraModTh]
construct a Thomason model structure
(modeling $(∞,0)$-categories)
on the category of small 2-categories and show it is Quillen equivalent to the Kan–Quillen model structure on simplicial sets.

Fiore–Paoli [\Thomncat] introduce a Thomason model structure on small strict $n$-fold categories
(obtained by iterating the internal category construction $n$ times starting from the category of sets)
and prove it is Quillen equivalent to simplicial sets.

Raptis [\Posets] constructs a Thomason model structure on the category of small posets,
shows it to be Quillen equivalent to the Thomason model structure on small categories,
and proves that the Thomason model structure is not cartesian.
Bruckner [\Acyclic] constructs a Thomason model structure on acyclic categories (categories without inverses and nonidentity endomorphisms).
Bruckner–Pegel [\ThCofibrant] give examples of cofibrant posets in the Thomason model structure.

Meier–Ozornova [\FPMC] show that the category of weak equivalences of a Barwick–Kan partial model category (Barwick–Kan [\Partial]) is Thomason-fibrant.
Meier [\FCFRC] shows that the category of weak equivalences of a fibration category is Thomason-fibrant.

\subsection Prerequisites

We assume familiarity with basics of the following topics from homotopy theory.
Appropriate references will be given throughout the text.

\li Simplicial homotopy theory, including simplicial sets, simplicial maps, simplicial weak equivalences,
Kan's subdivision–extension adjunction $\Sd⊣\Ex$,
the fundamental category–nerve adjunction $\fc⊣\N$.
See Gabriel–Zisman [\CFHT], Goerss–Jardine [\SHT].
\li Model categories, including model structures, Quillen adjunctions, transferred model structures, monoidal model categories.
See Quillen [\HoAlg], Hovey [\MC], Hirschhorn [\MCL], Barwick [\LR], Schwede–Shipley [\AMMMC].
\li Enriched categories, enriched operads, and algebras over enriched operads.
See Kelly [\EnrCat, \EnrOp].

\subsection Acknowledgments

I thank Daniel Grady for discussions connected to this paper.
I thank Lennart Meier and Viktoriya Ozornova for pointing out the results of Raptis [\Posets].
I thank Georges Maltsiniotis for ^!{couniversal weak equivalence}.
I thank Georges Maltsiniotis and Viktoriya Ozornova for providing feedback on a previous version of ^!{enriched Thomason Quillen equivalence}.
I thank the anonymous referee of the Journal of Pure and Applied Algebra for numerous suggestions and improvements.

\section The Thomason model structure and its properties

In this section we establish the properties of the Thomason model structure that we need later.

Recall the Kan–Quillen model structure on the category $\sSet$ of simplicial sets (Quillen [\HoAlg, Theorem~II.3.3]).
Recall the adjunction between Kan's subdivision and extension functors $$\Sd⊣\Ex,\qquad \Sd:\sSet→\sSet,\qquad \Ex:\sSet→\sSet \qquad \hbox{(Kan [\AHiii, §3, §6])}$$
and the adjunction between the fundamental category and nerve functors $$\fc⊣\N,\qquad \fc:\sSet→\Cat,\qquad \N:\Cat→\sSet \qquad \hbox{(Gabriel–Zisman [\CFHT, §II.4]),}$$
where $\Cat$ denotes the category of small categories and functors.

\proclaim Remark.
^^={properties of compact objects}
A simplicial set is a compact object in $\sSet$ if and only if it has finitely many nondegenerate simplices.
A small category is a compact object in $\Cat$ if and only if it admits a presentation with finitely many generators and relations.
Thus, the functors $\Sd$ and $\fc$ send representable simplicial sets $Δ^n$ to compact objects.
Therefore, $\Sd$ and $\fc$ preserve compact objects and $\Ex$ and $\N$ preserve filtered colimits.

Recall the definition of a transferred model structure.

\proclaim Definition.
(Crans [\QCMSS, Theorem~3.3], Hirschhorn [\MCL, Theorem~11.3.2].)
Suppose $C$ is a model category and $U:D→C$ is a right adjoint functor.
The ^={transferred model structure[|s]} on~$D$ (if it exists)
is the unique model structure whose weak equivalences and fibrations are created by the functor~$U$,
meaning a morphism~$f$ in~$D$ is a weak equivalence if and only if $U(f)$ is a weak equivalence in~$C$ and likewise for fibrations.

\proclaim Definition.
The ^={Thomason model structure} (Thomason [\CatCMC, Theorem~4.9], Cisinski [\Dwyer]) on the category $\Cat$ of small categories
is transferred (^!{transferred model structure}) along the right adjoint functor
$$\Ex^2 \N: \Cat→\sSet,$$
meaning its weak equivalences (^={Thomason weak equivalence[s|]}), fibrations, and acyclic fibrations are created by the functor $\Ex^2 \N$,
whereas cofibrations and acyclic cofibrations are defined using the left lifting property.

\proclaim Remark.
There are other model structures on $\Cat$, such as the model structure of Joyal–Tierney [\SSCS, Theorem~4],
in which weak equivalences are precisely equivalences of categories.
In this paper, $\Cat$ is always considered with the Thomason model structure.

\proclaim Proposition.
^^={properties of Thomason weak equivalences}
Weak equivalences in $\Cat$ are closed under compositions, transfinite compositions, filtered colimits, finite products, and satisfy the 2-out-of-3 property.

\proof Proof.
This follows from the analogous properties for $\sSet$
combined with the fact that the functor $\Ex^2\N$ preserves filtered colimits
because its left adjoint functor $\fc\Sd^2$ sends every simplex $Δ^m$ to a compact object in $\Cat$.

Recall the notion of an h-cofibration in a model category.
(The letter “h” stands for Hurewicz, meaning h-cofibrations abstract away some of the properties of Hurewicz cofibrations.)

\proclaim Definition.
(Grothendieck; Batanin–Berger [\HTAPM, Definition~1.1].)
A morphism $f:X→Y$ in a model category~$C$
(more generally, a relative category, i.e., a category equipped with a subcategory of weak equivalences)
is an ^={h-cofibration[|s]} if the cobase change functor along~$f$
$$f_!: X/C→Y/C$$
preserves weak equivalences, where $X/C$ denotes the undercategory (alias coslice category) of $X∈C$.
An ^={acyclic h-cofibration[|s]} is an h-cofibration that is also a weak equivalence.

Recall the following definition due to Cisinski, following Thomason.

\proclaim Definition.
(Cisinski [\Dwyer, Définition~1].)
A ^={Cisinski–Dwyer map[|s]} is an inclusion of a full subcategory
$i:A→B$
such that
\li The functor~$i$ is a sieve: if $Y∈A$ and $f:X→Y$ is a morphism in~$B$, then $X∈A$ and $f$ is a morphism in~$A$.
\li The inclusion $j:A→Z$ admits a retraction $r:Z→A$.
Here $Z⊂B$ is the cosieve generated by~$A$ in~$B$, i.e., the full subcategory of~$B$ comprising objects $Y∈B$ such that there is a morphism $X→Y$ in~$B$ for which $X∈A$.
\li The retraction~$r$ admits a natural transformation $ε:ir→\id_Z$.
\li The natural transformation $ε∘i$ is the identity natural transformation of functors $A→Z$.
\endlist
An ^={acyclic Cisinski–Dwyer map[|s]} is a Cisinski–Dwyer map that is also a weak equivalence in $\Cat$.

\proclaim Remark.
The original definition of Dwyer maps (Thomason [\CatCMC, Definition~4.1])
further requires that $i$ is left adjoint to~$r$, in which case $ε$ can be taken to be the counit.
This condition is too strong for Thomason [\CatCMC, Lemma~5.3.3] (whose proof is omitted there), which claims that Dwyer maps are closed under retracts.
Cisinski [\Dwyer, Lemme~1 and the two following paragraphs] constructs an example of a retract of a Dwyer map that is not a Dwyer map.
Cisinski [\Dwyer, Définition~1] defined Cisinski–Dwyer maps
and showed that Thomason's necessary criterion for cofibrations in $\Cat$ [\CatCMC, Proposition~5.4] and the proof of left properness of $\Cat$ [\CatCMC, Corollary~5.5]
are valid when Cisinski–Dwyer maps are substituted for Dwyer maps.

\proclaim Proposition.
^^={properties of Cisinski–Dwyer maps}
^{Cisinski–Dwyer maps} and ^{acyclic Cisinski–Dwyer maps} are closed under cobase changes, composition, transfinite composition, retracts, and products with a fixed object in $\Cat$.
All cofibrations in $\Cat$ are ^{Cisinski–Dwyer maps} (^!{Cisinski–Dwyer map}) and all ^{Cisinski–Dwyer maps} are ^{h-cofibrations} (^!{h-cofibration}) in $\Cat$.

\proof Proof.
^{Cisinski–Dwyer maps} are closed under retracts by Cisinski [\Dwyer, Lemme~4].
^{Cisinski–Dwyer maps} are closed under compositions and transfinite compositions by Cisinski [\Dwyer, Remarque] and Thomason [\CatCMC, Lemma~5.3].
^{Cisinski–Dwyer maps} are closed under cobase changes by Cisinski [\Dwyer, Remarque] and Thomason [\CatCMC, Proposition~4.3].
To show that ^{Cisinski–Dwyer maps} are closed under products with a fixed object $W∈\Cat$,
observe that if
$$(i:A→B,r:Z→A,ε:ir→\id_Z)$$
is a tuple exhibiting $i$ as a ^{Cisinski–Dwyer map}, then
$$(W⨯i:W⨯A→W⨯B,r:W⨯Z→W⨯A,W⨯ε:W⨯ir→W⨯\id_Z)$$
is a tuple exhibiting $W⨯i$ as a ^{Cisinski–Dwyer map}.
\ppar
All cofibrations in $\Cat$ are ^{Cisinski–Dwyer maps} by Cisinski [\Dwyer, Proposition~2].
To show that every ^{Cisinski–Dwyer map} $i:A→B$ is an ^{h-cofibration},
pick any weak equivalence $f:X→Y$,
together with a map $A→X$.
The induced map
$$X⊔_A B→Y⊔_A B$$
is a weak equivalence in $\Cat$, meaning
$$\N(X⊔_A B)→\N(Y⊔_A B)$$
is a weak equivalence in $\sSet$.
By Cisinski [\Dwyer, Remarque] and Thomason [\CatCMC, Proposition~4.3],
the latter morphism is weakly equivalent to
$$\N X⊔_{\N A}\N B→\N Y⊔_{\N A}\N B.$$
This morphism is a weak equivalence because $\N i:\N A→\N B$ is a cofibration of simplicial sets
and the model category $\sSet$ is left proper.
\ppar
Finally, ^{acyclic Cisinski–Dwyer maps} satisfy the same set of properties
by ^!{properties of Thomason weak equivalences} and the previously established properties of ^{Cisinski–Dwyer maps}.

\proclaim Definition.
(Barwick [\LR, Definition~1.3 (arXiv); 1.21 (journal)].)
A model category is ^={tractable}
if it is combinatorial (as defined by Jeffrey~H.~Smith, see Dugger [\Replacing, Definition~2.4])
and it admits a set of generating (acyclic) cofibrations with cofibrant domains.

\proclaim Remark.
^^={tractable cofibrations}
By Barwick [\LR, Corollary~1.12 (arXiv); 2.7 (journal)], a combinatorial model category that admits a set of generating cofibrations with cofibrant domains
is ^!{tractable}, i.e., admits a set of generating acyclic cofibrations with cofibrant domains.

The following definition is one of many ways to formalize the idea of a compactly generated model category.
The specific definition chosen here is motivated by the fact that it is well behaved with respect to transfers
of model structures, as shown in Pavlov–Scholbach [\Sym, Proposition~5.3(ii)].

\proclaim Definition.
(Pavlov–Scholbach [\Sym, Definition~2.1].)
A model category $C$ is ^={pretty small} if there is a cofibrantly generated model category structure~$D$ on~$C$
with the same weak equivalences, possibly smaller class of cofibrations,
and such that the domains and codomains of some set of generating cofibrations of~$D$ are compact objects.

Although the cited work of Thomason and Cisinski does prove the existence of the Thomason model structure, we find it beneficial to give a short modern proof,
while reusing only a fraction of their results.

Recall the transfer theorem for model structures, which we state in the special case of locally presentable categories.

\proclaim Proposition.
^^={existence of transferred model structures}
(Crans [\QCMSS, Theorem~3.3], Hirschhorn [\MCL, Theorem~11.3.2]; see also Garner–Kędziorek–Riehl [\LAMS] for a generalization.)
Suppose $U:D→C$ is a right adjoint functor between locally presentable categories
and $C$ is equipped with a combinatorial model structure with $I$ and $J$ as sets of generating (acyclic) cofibrations.
Then the ^{transferred model structure} (^!{transferred model structure}) on~$D$ exists if and only if
the functor~$U$ sends to weak equivalences all transfinite compositions of cobase changes of elements of $F(J)$,
where $F$ is the left adjoint of~$U$.
In this case, the model structure on~$D$ is combinatorial with $F(I)$ and $F(J)$ are sets of generating (acyclic) cofibrations.

\proclaim Theorem.
^^={existence of the Thomason model structure}
(Thomason [\CatCMC, Theorem~4.9], Cisinski [\Dwyer].)
The ^{Thomason model structure} exists and is a proper ^{tractable} (^!{tractable}) ^{pretty small} (^!{pretty small}) combinatorial model structure.

\proof Proof.
The functor $\Ex^2\N$ is a right adjoint functor between locally presentable categories that preserves filtered colimits and its left adjoint preserves compact objects.
Thus, by the transfer theorem (^!{existence of transferred model structures}),
it suffices to show that transfinite compositions of cobase changes of generating acyclic cofibrations
are weak equivalences in $\Cat$.
By ^!{properties of Cisinski–Dwyer maps}, generating acyclic cofibrations in $\Cat$ are ^{acyclic Cisinski–Dwyer maps}, which in their turn are closed under cobase changes and transfinite compositions.
Thus, the transferred model structure on $\Cat$ exists and is combinatorial.
By ^!{properties of Cisinski–Dwyer maps}, every cofibration in $\Cat$ is an ^{h-cofibration}, hence the model structure on $\Cat$ is left proper.
Since the model structure on $\sSet$ is right proper, so is the model structure on $\Cat$.
The model structure on $\Cat$ is ^{tractable} because the domains $\fc\Sd^2 ∂Δ^n$ and $\fc\Sd^2 Λ^n_k$ of generating (acyclic) cofibrations are cofibrant
since $\fc\Sd^2$ is a left Quillen functor and all simplicial sets are cofibrant.
The model structure on $\Cat$ is ^{pretty small} because the domains and codomains $\fc\Sd^2 Δ^n$ of generating (acyclic) cofibrations are compact
since $\fc\Sd^2$ preserves compact objects by ^!{properties of compact objects}
and the simplicial sets $∂Δ^n$, $Λ^n_k$, and $Δ^n$ have finitely many nondegenerate simplices, hence are compact.

\section Monoidal properties of the Thomason model structure

In this section, we investigate properties of the Thomason model structure connected to its cartesian monoidal structure.
An example due to Raptis [\Posets, §3] shows that the Thomason model structure is not cartesian.
Nevertheless, we can establish some other properties, which suffice to establish a model structure on modules over monoids in the Thomason model structure
(^!{Thomason modules}).

Recall the notion of a monoidal model category from Schwede–Shipley [\AMMMC].

\proclaim Definition.
(Schwede–Shipley [\AMMMC, Definition~2.1 (arXiv); Definition~3.1 (journal)].)
Suppose $M$ is a closed monoidal category equipped with a model structure.
We say that $M$ is a ^={monoidal model categor[y|ies]}
^^={monoidal model structure[|s]}
if the monoidal product functor $$M⨯M→M$$ is a left Quillen bifunctor.
The latter condition is also known as the ^={pushout product axiom[|s]}.
If the monoidal structure is cartesian, we talk about a ^={cartesian} model structure.

Often, an additional condition on the monoidal unit is included in the definition of a monoidal model category,
such as the unit axiom (Hovey [\MC, Lemma~4.2.7]), the strong unit axiom (Muro [\HTNOii, Definition~A.9]), or the very strong unit axiom (Batanin–Berger [\HTAPM, §1.20]).
In our case, the monoidal unit will always be cofibrant, which is the strongest of the unit axioms.

The following proposition is due to Raptis [\Posets, §3].

\proclaim Proposition.
(Raptis [\Posets, §3].)
The ^{Thomason model structure} is not ^{cartesian}.
More precisely, the pushout product of the acyclic cofibration $\{0\}→\{0→1\}$ with itself
is not a ^{Cisinski–Dwyer map} and therefore not a cofibration.

Below, ^!{flat} (flat model structure), ^!{h-monoidal} (h-monoidal model structure), ^!{monoid axiom} (the monoid axiom)
are stated for model categories equipped with a monoidal structure that need not satisfy the ^{pushout product axiom}.

\proclaim Definition.
(Pavlov–Scholbach [\Sym, Definition~3.2.4].)
A model category equipped with a monoidal structure is ^={flat}
if the pushout product of a cofibration and a weak equivalence
is a weak equivalence.

\proclaim Proposition.
^^={Thomason model structure is flat}
The ^{Thomason model structure} on $\Cat$ is ^{flat} (^!{flat}).

\proof Proof.
Suppose $A→B$ is a cofibration and $X→Y$ is a weak equivalence in $\Cat$.
Denote by $P→B⨯Y$ their pushout product:
$$\diagram[2em]A:0,0=$A⨯X$,B:3,0=$A⨯Y$,C:0,-3=$B⨯X$,D:3,-3=$B⨯Y$,E:2,-2=$P$,AC.lft=,BD.rt=,AB.top=,CD.bot=,BE.ulft=,CE.ulft=,ED.llft=;$$
The map $A⨯X→B⨯X$ is a ^{Cisinski–Dwyer map} and an ^{h-cofibration} by ^!{properties of Cisinski–Dwyer maps}.
The map $A⨯X→A⨯Y$ is a weak equivalence by ^!{properties of Thomason weak equivalences}.
Thus, the canonical map $B⨯X→P$ is a weak equivalence.
Since $B⨯X→B⨯Y$ is a weak equivalence by ^!{properties of Thomason weak equivalences},
by the 2-out-of-3 property the map $P→B⨯Y$ is also a weak equivalence.

\proclaim Corollary.
^^={acyclic pushout product}
The ^{Thomason model structure} on $\Cat$ (^!{Thomason model structure})
satisfies the acyclic part of the ^{pushout product axiom}:
the pushout product of a cofibration and an acyclic cofibration is a weak equivalence.

\proclaim Definition.
(Batanin–Berger [\HTAPM, Definition~1.11].)
A model category~$C$ equipped with a monoidal structure is ^={h-monoidal}
^^={h-monoidality}
if the monoidal product of any object $A∈C$
with an (acyclic) cofibration in~$C$
is an (acyclic) ^{h-cofibration} (^!{h-cofibration}).

Recall the Schwede–Shipley monoid axiom [\AMMMC, Definition~2.2 (arXiv); Definition~3.2 (journal)].

\proclaim Definition.
A model category~$C$ with a monoidal structure satisfies the ^={monoid axiom}
if every transfinite composition of cobase changes of tensor products of an object in~$C$ and an acyclic cofibration
is a weak equivalence in~$C$.

\proclaim Proposition.
^^={properties of the Thomason model structure}
The ^{Thomason model structure} on $\Cat$
(^!{Thomason model structure})
is a proper combinatorial
model category
that is ^{tractable} (^!{tractable}), ^{pretty small} (^!{pretty small}), ^{h-monoidal} (^!{h-monoidal}), ^{flat} (^!{flat}),
and satisfies the ^{monoid axiom} (^!{monoid axiom}).

\proof Proof.
By ^!{existence of the Thomason model structure}, the ^{Thomason model structure} is a proper combinatorial ^{tractable} ^{pretty small} model category.
\ppar
By ^!{properties of Cisinski–Dwyer maps},
cofibrations in $\Cat$ are ^{Cisinski–Dwyer maps},
the product of an object in $\Cat$ and a ^{Cisinski–Dwyer map} is again a ^{Cisinski–Dwyer map},
and every ^{Cisinski–Dwyer map} is an ^{h-cofibration},
hence the ^{Thomason model structure} is ^{h-monoidal}.
\ppar
The ^{monoid axiom} (^!{monoid axiom}) is established in Pavlov–Scholbach [\Sym, Lemma~3.2.3] (which does not require the ^{pushout product axioms}),
using the fact that the ^{Thomason model structure} is ^{h-monoidal} and ^{pretty small}.
We reproduce the proof of the cited result here:
by ^!{properties of Cisinski–Dwyer maps}, maps of the form $M⨯j$ (where $M∈\Cat$ and $j$ is an acyclic cofibration in $\Cat$) are ^{acyclic Cisinski–Dwyer maps},
and again by ^!{properties of Cisinski–Dwyer maps}, acyclic Cisinski–Dwyer maps are closed under cobase changes and transfinite compositions.

\proclaim Proposition.
^^={properties of the Thomason adjunction}
The Quillen equivalence
$$\fc\Sd^2⊣\Ex^2\N,\qquad \fc\Sd^2:\sSet→\Cat,\qquad \Ex^2\N:\Cat→\sSet$$
satisfies the conditions for a weak monoidal Quillen equivalence in the sense of Schwede–Shipley [\EMMC, Definition~3.6],
except for the nonacyclic part of the ^{pushout product axiom} for $\Cat$.
That is, $\sSet$ is a cartesian model category, $\Cat$ satisfies the acyclic pushout product axiom (^!{acyclic pushout product}),
the right adjoint $\Ex^2\N$ is a lax (in fact, strong) monoidal functor (i.e., preserves finite products),
the left adjoint $\fc\Sd^2$ preserves the monoidal unit (i.e., the terminal object),
and the canonical map
$$\fc\Sd^2(A⨯B)→\fc\Sd^2A⨯\fc\Sd^2B$$
is a weak equivalence for any simplicial sets $A$ and $B$.
Both functors in the adjunction preserve and reflect weak equivalences.
In particular, the unit and counit maps are weak equivalences.

\proof Proof.
The only statement that remains to be proved is that the comonoidal map
$$\fc\Sd^2(A⨯B)→\fc\Sd^2A⨯\fc\Sd^2B$$
is a weak equivalence for any simplicial sets $A$ and $B$.
This follows from the fact that the functor $\fc$ preserves finite products
and there is a natural weak equivalence $\Sd^2→\id$ induced by the last vertex map.

We conclude this section by constructing a model structure on the category of modules over a strict monoidal category, i.e., a monoid in the ^{Thomason model structure}.
Recall that any monoidal category can be strictified to a strict monoidal category.

\proclaim Proposition.
^^={Thomason modules}
Suppose $M$ is a monoid in the ^{Thomason model structure}, i.e., a strict monoidal category.
The category of (strict) modules over~$M$ admits a model structure transferred along the forgetful functor to the ^{Thomason model structure} (^!{Thomason model structure}).

\proof Proof.
(Compare Schwede–Shipley [\AMMMC, Remark~3.2 (arXiv); Remark 4.2 (journal)].)
By the transfer theorem (^!{existence of transferred model structures})
it suffices to show that the forgetful functor sends transfinite compositions of cobase changes of maps of the form $M⨯j$,
where $j∈J$ is a generating acyclic cofibration for the ^{Thomason model structure},
to weak equivalences in the ^{Thomason model structure}.
This is the content of the ^{monoid axiom} established in ^!{properties of the Thomason model structure}.

\section Nonexistence of cartesian model structures

In this section we show that any choice of a model structure on the category of small categories equipped with ^{Thomason weak equivalences}
must be a noncartesian model structure.
This means that the existing approaches to constructing model structures on monoids (Schwede–Shipley [\AMMMC]), enriched categories (Muro [\HTEC]),
or algebras over colored symmetric operads (Pavlov–Scholbach [\Adm])
do not work, since these make heavy use of the ^{pushout product axiom}.
Nevertheless, in the next section we do show the existence of model structures on monoids and enriched categories
by replicating the original argument of Thomason in this setting.

\proclaim Proposition.
^^={Thomason cofibrations}
Suppose the category of small categories is equipped with a model structure whose weak equivalences are ^{Thomason weak equivalences}.
Then the terminal category is cofibrant and at least one of the inclusions $\{0\}→\{0→1\}$, $\{1\}→\{0→1\}$ is an acyclic cofibration.

\proof Proof.
Factor the map $g:∅→\{0\}$
as a cofibration $p:∅→A$ followed by an acyclic fibration $q:A→\{0\}$.
Since $q$ is a ^{Thomason weak equivalence}, there is an object $a∈A$, which induces a morphism $ι:\{0\}→A$.
The maps $ι$ and $q$ exhibit $\{0\}$ as a retract of the cofibrant object~$A$,
hence the category $\{0\}$ is cofibrant.
In particular, acyclic fibrations must be surjective on objects.
\ppar
Factor the weak equivalence $\{0\}→\{0→1\}$ as an acyclic cofibration $p:\{0\}→F$ followed by an acyclic fibration $q:F→\{0→1\}$.
Any acyclic fibration has the right lifting property with respect to the cofibration $∅→\{0\}$, hence $q$ is surjective on objects.
Since $q$ is a ^{Thomason weak equivalence}, the arrow $0→1$ must be in the image of~$q$.
\ppar
Introduce a preorder relation~$P$ on the set~$S$ of objects of~$F$: we have $x≤y$ if there is an arrow $x→y$.
Consider also an equivalence relation~$R$ on the same set~$S$: we have $xRy$ if $x≤y$ and $y≤x$.
The preorder~$P$ becomes an order on the quotient $S/R$.
By Szpilrajn's extension theorem [\EOP], we can extend the partial order on $S/R$ to a total order on~$S/R$,
which induces a total preorder~$Q$ on~$S$.
\ppar
The equivalence class~$E$ (with respect to the equivalence relation~$R$ on~$S$)
of~$p(0)∈F$ must necessarily have at least incoming or outgoing arrow connecting it (and hence connecting~$p(0)$) to a different equivalence class~$E'$,
since otherwise $E$ is a connected component of~$F$, which is impossible because such a connected component cannot contain an object $z∈F$ such that $q(z)=1$.
Since $Q$ is a total order, we have $E<E'$ or $E'<E$, so we analyze each of those cases separately.
\ppar
Suppose $z→p(0)$ is an arrow in~$F$ such that $z<p(0)$ with respect to~$Q$.
Consider the inclusion $\{0→1\}→F$ that sends $0→1$ to $z→p(0)$
and the retraction $F→\{0→1\}$ that sends $z↦0$ if $z<p(0)$ and $z↦1$ if $z≥p(0)$.
The definition of $Q$ guarantees that the latter construction yields a functor.
We exhibited the injective functor $\{0\}→\{0→1\}$ that sends $0↦1$ as a codomain retract of the inclusion $\{0\}→F$.
Thus, the inclusion $\{1\}→\{0→1\}$ is a cofibration.
\ppar
Suppose $p(0)→z$ is an arrow in~$F$ such that $p(0)<z$ with respect to~$Q$.
Consider the inclusion $\{0→1\}→F$ that sends $0→1$ to $p(0)→z$
and the retraction $F→\{0→1\}$ that sends $z↦0$ if $z≤p(0)$ and $z↦1$ if $z>p(0)$.
The definition of $Q$ guarantees that this is a functor.
We exhibited the inclusion $\{0\}→\{0→1\}$ as a codomain retract of the inclusion $\{0\}→F$.
\ppar
Thus, at least one of the inclusions $\{0\}→\{0→1\}$, $\{1\}→\{0→1\}$ is a cofibration.

I thank Georges Maltsiniotis for suggesting the choice of the map $A→D$ in the following proposition.

\proclaim Proposition.
^^={couniversal weak equivalence}
(Maltsiniotis, 2022.)
Consider the posets $B=\{0→1\}^3$, $A=B∖\{111\}$, $C=A∖\{000\}$.
Turn these posets into categories,
and consider also the categories given by pushouts in the category of small categories $D=A⊔_C A$ and $E=D⊔_A B$.
The canonical map $D→E$ is not a ^{Thomason weak equivalence}:
$$\sqcd{
C&\mapright{}&A\cr
\mapdown{}&&\mapdown{}\cr
A&\mapright{}&D\cr
\mapdown{}&&\mapdown{}\cr
B&\mapright{}&E.\cr
}$$

\proof Proof.
The inclusion~$ι$ of posets into categories preserves the pushout $A⊔_C A$,
i.e., the canonical map $$ιA⊔_{ιC}ιA→ι(A⊔_C A)$$ is an isomorphism because $C⊂A$ is an upward closed subset.
Henceforth, we omit~$ι$ from the notation.
Thus, the category~$D$ is the nerve of the poset whose underlying set is $A⊔\{000'\}$, where $000'$
is ordered in the same way as $000$, and the elements $000$ and $000'$ are incomparable.
Furthermore, the nerve functor preserves the pushout $A⊔_C A$,
i.e., the canonical map $$\N A⊔_{\N C}\N A→\N(A⊔_C A)$$
is an isomorphism.
Thus, the homotopy type of $D$ can be computed as the homotopy pushout of $\N A←\N C→\N A$.
The nerve of~$A$ is contractible.
The nerve of~$C$ is six 1-simplices connected together, i.e., a circle.
Thus, $\N D$ is weakly equivalent to the 2-sphere.
\ppar
Consider now the pushout $E=D⊔_A B$.
The category $E$ has the same objects as~$B$, together with the additional object $000'$ coming from~$D$.
We claim that the object $111∈E$ is a terminal object.
Suppose $f,g:a→111$ are parallel morphisms from some object $a∈E$.
If $a∈B$, then $f$ and $g$ necessarily come from~$B$,
where they are equal.
The only other option is $a=000'$.
In this case, there are six potential morphisms $000'→111$
given by composing one of the six edges $b→c$ in~$C$
with the unique morphisms $000'→b$ and $c→111$.
All six morphisms coincide, thanks to relations
like $000'→001→011$ being equal to $000'→010→011$
and $001→011→111$ being equal to $001→101→111$.
Thus, $111∈E$ is a terminal object and $E$ is weakly contractible.

We remark that working in the pushout of categories $D ⊔_A B$ is crucial for the above argument.
In particular, if we worked with the pushout of nerves $\N D ⊔_{\N A} \N B$,
then the above argument would construct simplicial homotopies between the corresponding triples of edges,
which glue together into a noncollapsible 2-sphere.
In the pushout of categories $D ⊔_A B$ this 2-sphere is collapsed to a single morphism $000'→111$.

\proclaim Proposition.
^^={nonexistence of cartesian structures}
Suppose the category of small categories is equipped with a model structure whose weak equivalences are ^{Thomason weak equivalences}.
Then this model structure is not cartesian.

\proof Proof.
By ^!{Thomason cofibrations}, at least one of the inclusions $\{0\}→\{0→1\}$, $\{1\}→\{0→1\}$ is a cofibration in this model structure.
The two cases are symmetric, so without loss of generality we assume $\{0\}→\{0→1\}$ is a cofibration.
The pushout product of three copies of this cofibration is the nerve of inclusion of posets $$ι:\{0→1\}^3∖\{111\}⊂\{0→1\}^3,$$
i.e., the map $ι:A→B$ from ^!{couniversal weak equivalence}.
If the ^{pushout product axiom} is satisfied, the map~$ι$ must be an acyclic cofibration, in particular, its cobase changes must be ^{Thomason weak equivalences}.
However, ^!{couniversal weak equivalence} constructs a cobase change of~$ι$ that is not a ^{Thomason weak equivalence}.

\section Main theorems

In this section we construct the enriched Thomason model structure on small 2-categories (^!{enriched Thomason model structure})
and prove it is Quillen equivalent to the Bergner model structure on small simplicial categories (^!{enriched Thomason Quillen equivalence}).
Thus, strict 2-categories equipped with the enriched Thomason model structure provide another model for $(∞,1)$-categories.
We start with the simpler case of a model structure on monoids in the Thomason model structure (^!{model structure on monoids}).
The established machinery for constructing model structures on monoids (Schwede–Shipley [\AMMMC])
is not applicable because ^!{nonexistence of cartesian structures} shows that any approach based on the ^{pushout product axiom} cannot work.
However, we do make use of the Schwede–Shipley filtration (^!{monoid filtration}) on cobase changes of free morphisms of monoids.

\proclaim Definition.
^^={monoid conventions}
Suppose $C$ and $D$ are locally presentable closed symmetric monoidal categories.
Denote by $$\Mon_C$$ the category of monoids in~$C$.
Denote by
$$\Free⊣\Under: C ⇄ \Mon_C$$
the free-forgetful adjunction between $C$ and monoids in~$C$.
We also write $\Free_C⊣\Under_C$ if there is more than one possibility for~$C$.
Given an adjunction $L⊣R:C⇄D$,
where the right adjoint functor~$R$ is strong monoidal
(in our case: $\fc\Sd^2⊣\Ex^2\N:\sSet⇄\Cat$), denote by
$$\FreeM⊣\UnderM: \Mon_C ⇄ \Mon_D$$
the adjunction between monoids in~$C$ and monoids in~$D$,
where the right adjoint $\UnderM$ applies the functor~$R$ to the underlying object of a monoid in~$D$.

\proclaim Definition.
Given a morphism $f:K→L$ in a cocomplete symmetric monoidal category~$C$
and $i≥0$, denote by $$f^{◻i}:f^{⊡i}→L^{⊗i}$$ the $i$-fold pushout product of~$f$ with itself.

\proclaim Remark.
The notation $f^{⊡i}$ refers to the domain of $f^{◻i}$.
If $C$ is the category of sets with its cartesian monoidal structure and $f$ is an inclusion of sets, then
$f^{⊡i}$ is the set $L^{⊗i}∖(L∖K)^{⊗i}$.

Recall the following presentation of cobase changes of free morphisms of monoids due to Schwede–Shipley [\AMMMC, proof of Lemma~5.2 (arXiv); 6.2 (journal)].

\proclaim Proposition.
^^={monoid filtration}
(Schwede–Shipley [\AMMMC, proof of Lemma~5.2 (arXiv); 6.2 (journal)].)
Assuming the notation of ^!{monoid conventions},
given a morphism $f:K→L$ in~$C$ and a morphism $\Free(K)→X$ in $\Mon_C$,
in the pushout square
$$\cd{
\Free(K)&\mapright{}&X\cr
\mapdown{}&&\mapdown{}\cr
\Free(L)&\mapright{}&P\cr
}$$
the map $\Under(X)→\Under(P)$ is the transfinite composition of the chain
$$\Under(X)=P_0→P_1→P_2→⋯→P_∞=\Under(P),$$
where for every $i>0$ the morphism $P_{i-1}→P_i$ fits into the pushout square
$$\cd{
f^{⊡i}⊗X^{⊗(i+1)}&\mapright{}&P_{i-1}\cr
\lmapdown{f^{◻i}⊗X^{⊗(i+1)}}&&\mapdown{}\cr
L^{⊗i}⊗X^{⊗(i+1)}&\mapright{}&P_i,\cr
}$$
where the attaching map on top is defined inductively by using the previously defined maps $$L^{⊗j}⊗X^{⊗(j+1)}→P_j→P_{i-1}$$ for $j<i$.

\proclaim Proposition.
^^={monoidal h-cofibrations}
Assuming the notation of ^!{monoid conventions},
given $n≥0$ and a commuting diagram of pushout squares in the category $\Mon_\Cat$
$$\cd{
\Free(\fc\Sd^2 ∂Δ^n)&\mapright{a}&X&\mapright{}&Y\cr
\mapdown{}&&\mapdown{}&&\mapdown{}\cr
\Free(\fc\Sd^2 Δ^n)&\mapright{}&P&\mapright{}&Q,\cr
}$$
if the map $X→Y$ is a ^{Thomason weak equivalence},
then so is the map $P→Q$.
That is to say, applying $\Free∘\fc∘\Sd^2$ to $∂Δ^n→Δ^n$ yields an ^{h-cofibration} (^!{h-cofibration}) in $\Mon_\Cat$.

\proof Proof.
Our proof follows the analogous proof of Thomason [\CatCMC, Proposition~4.3], combining it with the filtration of ^!{monoid filtration}.
Recall that the inclusion $\Sd^2 ∂Δ^n→\Sd^2 Δ^n$ is the nerve of a map of posets $ι:A→B$,
which factors as $$\fc\Sd^2 ∂Δ^n=A\lto5{κ}W\lto5{μ}B=\fc\Sd^2 Δ^n,$$ where $W$ is the cosieve generated by~$A$ in~$B$.
Recall also that $W$ exhibits $A$ as a categorical analogue of a neighborhood deformation retract;
indeed, the inclusion $A→W$ is isomorphic to the inclusion $A⨯\{0\}→A⨯\{0→1\}$,
so in particular, it admits a retraction $W→A$ and the composition $W→A→W$ admits a natural transformation to the identity functor.
\ppar
Consider the induced diagram of pushout squares:
$$\cd{
\Free A&\mapright{}&X&\mapright{}&Y\cr
\lmapdown{\Free κ}&&\mapdown{}&&\mapdown{}\cr
\Free W&\mapright{}&R&\mapright{}&S\cr
\lmapdown{\Free μ}&&\mapdown{}&&\mapdown{}\cr
\Free B&\mapright{}&P&\mapright{}&Q.\cr
}$$
The rest of the proof proceeds in two steps: first, we show that
the map $R→S$ is a ^{Thomason weak equivalence}, and then we prove the same claim for the map $P→Q$.
\ppar
To show that $R→S$ is a ^{Thomason weak equivalence}, we exhibit the map $\Free A→\Free W$
as a monoidal analogue of a ^{Cisinski–Dwyer map}.
First, the retraction $W→A$ induces a retraction $\Free W → \Free A$ and therefore a retraction $R→X$.
In particular, the composition $X→R→X$ is the identity map.
To show that the composition $R→X→R$ is a ^{Thomason weak equivalence},
we apply the functor~$\Under$
and construct a natural transformation
$$ε:\Under(R→X→R)→\id_{\Under R}$$ using the filtration of ^!{monoid filtration}.
Thus, present the functor $\Under X→\Under R$ as the transfinite composition
$$\Under X=R_0→R_1→R_2→⋯→R_∞=\Under R,$$
where for every $i>0$ the map $R_{i-1}→R_i$ is a cobase change of the map $$X^{⨯(i+1)}⨯(κ^{⊡i}→W^{⨯i}),$$ where $κ:A→W$ is the inclusion map.
The natural transformation
$$ε:I⨯R_∞→R_∞,\qquad I=\{0→1\}$$
is constructed as the colimit of natural transformations
$$ε_i:I⨯R_i→R_i,$$
which are constructed by induction on~$i≥0$.
We take $ε_0:I⨯X→X$ to be the projection functor.
This also ensures that $ε∘χ$ is the identity natural transformation of functors $X→R$, where $χ:X→R$ is the inclusion map,
in analogy to the definition of a ^{Cisinski–Dwyer map}.
By inductive assumption, we have already constructed the natural transformation $$ε_{i-1}:I⨯R_{i-1}→R_{i-1}.$$
The natural transformation $ε_i$ is constructed as follows:
$$\eqalign{
ε_i:I⨯R_i&=I⨯(R_{i-1}⊔_{X^{⨯(i+1)}⨯κ^{⊡i}}X^{⨯(i+1)}⨯W^{⨯i})\cr
&{}≅(I⨯R_{i-1})⊔_{I⨯X^{⨯(i+1)}⨯κ^{⊡i}}(I⨯X^{⨯(i+1)}⨯W^{⨯i})\cr
&{}→R_{i-1}⊔_{X^{⨯(i+1)}⨯κ^{⊡i}}X^{⨯(i+1)}⨯W^{⨯i}=R_i,\cr
}$$
where the last map has as its first component~$ε_{i-1}$
and the other two components are induced by the identity on $X^{⨯(i+1)}$ and the natural transformation of functors
$$α_i:(W^{⨯i}→A^{⨯i}→W^{⨯i})⇒\id_{W^{⨯i}}$$
given by taking the product of the identity functor on~$A^{⨯i}$
with the natural transformation
$$(I^{⨯i}→\{0\}^{⨯i}→I^{⨯i})⇒\id_{I^{⨯i}},\qquad x↦(0^{⨯i}→x).$$
By inspection, the map $α_i$ is compatible with the map $ε_{i-1}$, so we indeed get a morphism of pushouts $ε_i:I⨯R_i→R_i$.
This proves that the map $X→R$ is a ^{Thomason weak equivalence}.
The same argument also proves that $Y→S$ is a weak equivalence.
By the 2-out-of-3 property, the morphism $R→S$ is also a weak equivalence.
\ppar
We now proceed to show that the map $P→Q$ is a ^{Thomason weak equivalence}.
Following Thomason, denote by $V⊂B$ the full subcategory given by the objects of~$B$ that are not in~$A$.
The inclusions $V→B$ and $W∩V→W$ are cosieves.
Here $V∖W$ has a single vertex, the barycenter of $\Sd^2 Δ^n$, which is the initial object of~$V$.
The inclusion $W∩V→W$ is isomorphic to the inclusion $$A⨯\{1\}→A⨯\{0→1\}≅W.$$
As observed by Thomason [\CatCMC, Proof of Proposition~4.3], the commutative square
$$\sqcd{
W∩V&\mapright{}&W\cr
\lmapdown{λ}&&\mapdown{μ}\cr
V&\mapright{}&B\cr
}$$
is cocartesian because the left and top maps are inclusions of cosieves.
Thus, we have a commutative diagram of cocartesian squares
$$\cd{
\hskip-.35in \Free(W∩V)&\mapright{}&\Free W&\mapright{}&R&\mapright{}&S\cr
\lmapdown{\Free λ}&&\lmapdown{\Free μ}&&\mapdown{}&&\mapdown{}\cr
\Free V&\mapright{}&\Free B&\mapright{}&P&\mapright{}&Q.\cr
}$$
From now on, we ignore the second column with the map $\Free μ$
and analyze the maps $R→P$ and $S→Q$ using the filtration of ^!{monoid filtration} with respect to the map $\Free λ$.
The $i$th step in the filtration of the map $\Under R→\Under P$ has the form
$$\cd{
R^{⨯(i+1)}⨯λ^{⊡i}&\mapright{}&P_{i-1}\cr
\lmapdown{R^{⨯(i+1)}⨯λ^{◻i}}&&\mapdown{}\cr
R^{⨯(i+1)}⨯V^{⨯i}&\mapright{}&P_i,\cr
}$$
where $λ:W∩V→V$ is the inclusion map.
The left map is a cosieve because $V$ adds a single vertex to $W∩V$ given by the barycenter,
so in the pushout product, $V^{⨯i}$ adds a single vertex to $λ^{⊡i}$ given by the barycenter in every factor,
and such a vertex does not admit any morphisms from $λ^{⊡i}$.
Next, we show that the inclusion $$Y=R^{⨯(i+1)}⨯λ^{⊡i}→P_{i-1}$$ is a cosieve.
Given an object $$z=r_0v_1r_1v_2⋯r_{i-1}v_ir_i∈Y$$ (where $r_k∈R$, $v_k∈V$, and juxtaposition denotes the multiplication operation),
morphisms of the form $z→z'$ in $P_{i-1}$ are by construction monoidal products of morphisms $$r_0→r'_0, v_1→v'_1, r_1→r'_1, …$$
because $V∩W$ is a cosieve in~$W$.
Thus, for such a morphism $z→z'$ the object~$z'$ is in the image of~$Y$ because $V∩W$ is a cosieve in~$V$.
By the same argument, the map $$Y=R^{⨯(i+1)}⨯λ^{⊡i}→Q_{i-1}$$ is also a cosieve inclusion.
\ppar
As observed by Thomason [\CatCMC, Proof of Proposition~4.3], a pushout square of cosieve inclusions in $\Cat$
is a homotopy pushout square in $\Cat$.
Thus, the cobase change square yielding $P_{i-1}→P_i$ is homotopy cocartesian.
The same argument shows that the square
$$\cd{
S^{⨯(i+1)}⨯λ^{⊡i}&\mapright{}&Q_{i-1}\cr
\lmapdown{S^{⨯(i+1)}⨯λ^{◻i}}&&\mapdown{}\cr
S^{⨯(i+1)}⨯V^{⨯i}&\mapright{}&Q_i,\cr
}$$
is homotopy cocartesian.
The morphism $R→S$ induces a natural transformation from the former square to the latter square.
In this natural transformation, the component $P_{i-1}→Q_{i-1}$ is a weak equivalence by induction,
with the base case $i=1$ yielding the map $R→S$, which is a ^{Thomason weak equivalence} by the first part of the proof.
The other two components are given by taking the product of the map $R^{⨯(i+1)}→S^{⨯(i+1)}$ with the category
$λ^{⊡i}$ respectively $V^{⨯i}$.
By ^!{properties of Thomason weak equivalences}, these products are weak equivalences because $R→S$ is a weak equivalence.
Therefore, all three components are weak equivalences.
The homotopy pushout of weak equivalences is again a weak equivalence.
Hence the map $P_i→Q_i$ is a weak equivalence in $\Cat$.
Therefore, the map $P→Q$ is also a weak equivalence in $\Cat$ by ^!{properties of Thomason weak equivalences}.

\proclaim Theorem.
^^={model structure on monoids}
The category $\Mon_\Cat$ of monoid objects in small categories (i.e., strict monoidal categories)
admits a proper model structure transferred
(^!{transferred model structure})
along the right adjoint forgetful functor $\Mon_\Cat→\Cat$
from the ^{Thomason model structure} (^!{Thomason model structure}) on $\Cat$.

\proof Proof.
The set of generating (acyclic) cofibrations in $\Mon_\Cat$
is given by $\Free(\fc\Sd^2(I))$ respectively $\Free(\fc\Sd^2(J))$.
Weak equivalences in $\Cat$ are closed under transfinite compositions
by ^!{properties of Thomason weak equivalences}.
Thus, by Pavlov–Scholbach [\Sym, Lemma~2.5(iv)]
the class of h-cofibrations in $\Mon_\Cat$ is a weakly saturated class,
i.e., is closed under cobase changes, transfinite compositions, and retracts.
Since $J$ is a subset of the weak saturation of~$I$,
we deduce that $\Free(\fc\Sd^2(J))$ is a subset of the weak saturation of $\Free(\fc\Sd^2(I))$;
in particular, it consists of h-cofibrations by ^!{monoidal h-cofibrations}.
This means that cobase changes of elements of $\Free(\fc\Sd^2(J))$ are weak equivalences.
Hence so are their transfinite compositions (by ^!{properties of Thomason weak equivalences}) and retracts.
By Hirschhorn [\MCL, Theorem~11.3.2] this proves the existence of the transferred model structure.
Since cofibrations are h-cofibrations, the transferred model structure is left proper,
whereas right properness is inherited from $\Cat$.
\ppar
An alternative proof can be given by citing Lurie [\HTT, Proposition~A.2.6.15],
which requires the class of weak equivalences in $\Mon_\Cat$ to be perfect (follows from ^!{properties of Thomason weak equivalences}),
the generating cofibrations to be h-cofibrations (holds by ^!{monoidal h-cofibrations}),
and morphisms with the right lifting property with respect to all generating cofibrations to be weak equivalences (holds by adjunction between $\Cat$ and $\Mon_\Cat$).

\proclaim Definition.
If $V$ is a monoidal category, then $\ECat_V$ denotes the category of small $V$-enriched categories and $V$-enriched functors.
In the special case of $V=\Cat$ with the cartesian monoidal structure, we talk about small ^={2-categor[ies|y]}.

\proclaim Remark.
The category $\ECat_V$ is locally presentable whenever $V$ is a locally presentable closed monoidal category.
For a detailed proof, see Kelly–Lack [\KL, Theorem~4.5].

\proclaim Definition.
^^={Dwyer–Kan maps}
(Muro [\HTEC].)
Suppose $V$ is a model category with a monoidal structure and $F:C→D$ is a $V$-enriched functor between $V$-enriched categories.
\li The functor~$F$ is {\it essentially surjective\/} if it becomes an essentially surjective functor in the usual sense
after applying the functor $A↦[1,A]$ to every hom-object, where $[-,-]$ denotes the hom-set in the homotopy category of~$V$ and $1$ denotes the monoidal unit in~$V$.
\li The functor~$F$ is a ^={Dwyer–Kan equivalence[|s]}
if the functor~$F$ is essentially surjective, and for every pair of objects $x,y∈C$ the induced morphism $C(x,y)→D(Fx,Fy)$ is a weak equivalence in~$V$.
\li The functor~$F$ is a ^={Dwyer–Kan acyclic fibration[|s]}
if the functor~$F$ is surjective on objects, and for every pair of objects $x,y∈C$ the induced morphism $C(x,y)→D(Fx,Fy)$ is an acyclic fibration in~$V$.

\proclaim Definition.
(Muro [\HTEC].)                                                                                                                                        
Suppose $V$ is a model category with a monoidal structure.
The ^={Dwyer–Kan model structure[|s]} on $\ECat_V$ (if it exists) is a model structure with weak equivalences and acyclic fibrations as in ^!{Dwyer–Kan maps}.

\proclaim Theorem.
The category $\ECat_\Cat$ of small ^{2-categories} admits a left proper combinatorial model structure
(the ^={enriched Thomason model structure})
whose weak equivalences are ^{Dwyer–Kan equivalences}
and acyclic fibrations are ^{Dwyer–Kan acyclic fibrations}
as in ^!{Dwyer–Kan model structure}, with $V=\Cat$ being equipped with the ^{Thomason model structure}.

\proof Proof.
Like in the proof of ^!{model structure on monoids}, we invoke Lurie [\HTT, Proposition~A.2.6.15].
As before, the class of Dwyer–Kan weak equivalences is a perfect class by ^!{properties of Thomason weak equivalences} and Muro [\HTEC, Proposition~9.2].
Take $$I'=\{∅→\{0\}\} ∪ T_{0,1}(I)$$ as a set of generating cofibrations for $\ECat_\Cat$.
Here $T_{0,1}:V→\ECat_V$ is a functor that sends an object $A∈V$ to the category $T_{0,1}(A)$ defined as follows.
The set of objects of~$T_{0,1}(A)$ is $\{0,1\}$.
The only nontrivial hom-object of $T_{0,1}(A)$ is $T_{0,1}(A)(0,1)=A$.
By Muro [\HTEC, Corollary~4.8], morphisms with the right lifting property with respect to~$I'$
are Dwyer–Kan equivalences.
Thus, it remains to show that elements of $I'$ are ^{h-cofibrations}.
The morphism $∅→\{0\}$ is an ^{h-cofibration} because Dwyer–Kan equivalences are closed under disjoint unions.
\ppar
To show that the elements of $T_{0,1}(I)$ are ^{h-cofibrations}, we use the same argument as in ^!{monoidal h-cofibrations}.
Given $n≥0$ and a commuting diagram of pushout squares in the category $\ECat_\Cat$
$$\cd{
T_{0,1}(\fc\Sd^2 ∂Δ^n)&\mapright{}&X&\mapright{}&Y\cr
\lmapdown{T_{0,1}(ι)}&&\mapdown{}&&\mapdown{}\cr
T_{0,1}(\fc\Sd^2 Δ^n)&\mapright{}&P&\mapright{s}&Q,\cr
}$$
we claim that if the map $X→Y$ is a Thomason weak equivalence,
then so is the map $s:P→Q$.
That is to say, applying $T_{0,1}∘\fc∘\Sd^2$ to $∂Δ^n→Δ^n$ yields an ^{h-cofibration} (^!{h-cofibration}) in $\ECat_\Cat$.
By Muro [\HTEC, Proposition~9.1], homotopically essentially surjective functors are closed under cobase changes.
Thus, it suffices to show that $P→Q$ is homotopically fully faithful,
i.e., for any objects $x,y∈P$, the induced map $$P(x,y)→Q(s(x),s(y))$$ is a ^{Thomason weak equivalence}.
(Recall that $X→P$ induces an identity map on the sets of objects.)
\ppar
Given a morphism $f:K→L$ in $\Cat$,
consider the following cocartesian square in $\ECat_\Cat$:
$$\cd{
T_{0,1}(K)&\mapright{}&X\cr
\lmapdown{T_{0,1}(f)}&&\mapdown{}\cr
T_{0,1}(L)&\mapright{}&P.\cr
}$$
Here the top attaching map sends $0↦a$, $1↦b$ for some objects $a,b∈X$.
As shown by Muro [\HTEC, §5], unfolding ^!{monoid filtration} yields a presentation of the morphism $X(x,y)→P(x,y)$ (for some $x,y∈X$)
as the transfinite composition of cobase changes of morphisms (Muro [\HTEC, (5.5)])
$$X(b,y)⨯X(b,a)^{⨯(i-1)}⨯X(x,a)⨯f^{◻i}.$$
\ppar
As in ^!{monoidal h-cofibrations}, we factor the map~$ι=\fc\Sd^2 δ_n$ as
$$A=\fc\Sd^2 ∂Δ^n \lto5{κ} W \lto5{μ} \fc\Sd^2 Δ^n=B$$
and apply the above filtration first to $f=κ$ and then to $f=λ$, where $λ:W∩V→V$ is defined in ^!{monoidal h-cofibrations}.
If we now compare the above map to the map $X^{⨯(i+1)}⨯κ^{◻i}$ (respectively $R^{⨯(i+1)}⨯λ^{◻i}$) used in the single-object case in ^!{monoidal h-cofibrations},
we see that the object $X^{⨯(i+1)}$ (respectively $R^{⨯(i+1)}$) is manipulated formally in the entire proof.
The remainder of the proof now proceeds identically to ^!{monoidal h-cofibrations},
replacing $X^{⨯(i+1)}$ with $$X(b,y)⨯X(b,a)^{⨯(i-1)}⨯X(x,a)$$ everywhere in the first part of the proof
and
replacing $R^{⨯(i+1)}$ with $$R(b,y)⨯R(b,a)^{⨯(i-1)}⨯R(x,a)$$ everywhere in the second part of the proof.

\proclaim Remark.
^^={simplicial categories are left proper}
Since ^!{enriched Thomason model structure} shows that $\ECat_\Cat$ is left proper,
one can ask a similar question about the Bergner model structure on $\ECat_\sSet$.
Indeed, Lurie [\HTT, Proposition~A.3.2.4] (available in arXiv v1 from 2006)
and Cisinski–Moerdijk [\SOpII, Corollary~8.10] show that the Bergner model structure is left proper.

\proclaim Remark.
One might wonder whether it may be possible to extend ^!{model structure on monoids} to construct
a model structure on algebras over operads.
In this case, the Schwede–Shipley filtration of ^!{monoid filtration}
is replaced by the Elmendorf–Mandell filtration [\FreePushouts, §12],
which replaces the map
$$f^{◻i}⊗X^{⊗(i+1)}:f^{⊡i}⊗X^{⊗(i+1)}→L^{⊗i}⊗X^{⊗(i+1)}$$
by the map
$$f^{◻i}⊗\Env(O,X)_i:f^{⊡i}⊗\Env(O,X)_i→L^{⊗i}⊗\Env(O,X)_i,$$
where $\Env(O,X)$ is the enveloping operad of the $O$-algebra~$X$, defined using a universal property,
and $\Env(O,X)_i$ denotes its object of operations of arity~$i$.
For symmetric operads, we must also mod out by the action of the symmetric group~$Σ_i$.
In the case of monoids (i.e., $O$ is the associative nonsymmetric operad), we have $\Env(O,X)=X^{⊗(i+1)}$,
which recovers the Schwede–Shipley filtration.
The proof of ^!{model structure on monoids} used the fact that if $X→Y$ is a weak equivalence,
then so is $\Env(O,X)→\Env(O,Y)$.
This is trivial if $\Env(O,X)=X^{⊗(i+1)}$,
but for an arbitrary operad~$O$ we have no control over $\Env(O,X)$
unless we impose very strong conditions on~$O$ such as cofibrancy.
It seems plausible that arguments of this type could be used to construct
an enriched Thomason model structure on {\it nonsymmetric\/} colored operads enriched in $\Cat$
along the lines of the model structure of Cisinski–Moerdijk [\SOpII, §1] for symmetric simplicial operads.
The symmetric case creates further problems, for example, by Hackney–Robertson–Yau [\Properness, §4],
the model category of single-colored symmetric simplicial operads is not left proper,
whereas our strategy relies on showing that cofibrations are h-cofibrations, which implies left properness.
The counterexample uses operations of arity~0 in an essential way,
which leaves open the question whether {\it reduced\/} colored symmetric operads enriched in $\Cat$ admit an enriched Thomason model structure.

\proclaim Theorem.
^^={enriched Thomason Quillen equivalence}
The Quillen equivalence of ^!{properties of the Thomason adjunction}
induces a Quillen equivalence
$$\FreeM⊣\UnderM,\qquad \FreeM:\ECat_\sSet→\ECat_\Cat,\qquad \UnderM:\ECat_\Cat→\ECat_\sSet,$$
where the right adjoint functor~$\UnderM=\ECat_{\Ex^2\N}$ applies the functor $\Ex^2\N$ to every hom-category in a given ^{2-category}.
That is to say, the model category of small simplicial categories (Bergner [\SCat])
and small ^{2-categories} with the ^{enriched Thomason model structure} (^!{enriched Thomason model structure})
are Quillen equivalent.
Restricting to categories with a single object yields a Quillen equivalence
$$\FreeM⊣\UnderM,\qquad \FreeM:\Mon_\sSet→\Mon_\Cat,\qquad \UnderM:\Mon_\Cat→\Mon_\sSet.$$

\proof Proof.
The functor~$\UnderM$ is well defined because the functor $\Ex^2\N$ preserves finite products up to a natural isomorphism.
The functor~$\UnderM$ preserves small limits and filtered colimits since these limits and colimits are computed on the level of underlying graphs by Muro [\HTEC, Proposition~2.6].
Thus, $\UnderM$ is an accessible continuous functor between locally presentable categories, hence a right adjoint functor by the adjoint functor theorem.
\ppar
By ^!{properties of the Thomason adjunction}, the functor $\UnderM$ preserves and reflects weak equivalences.
The explicit description of generating cofibrations in $\ECat_\sSet$
(Bergner [\SCat, Proposition~3.2])
as maps $∅→\{0\}$ and $T_{0,1}(∂Δ^n→Δ^n)$
immediately implies that $\FreeM$ preserves cofibrations.
The explicit description of generating acyclic cofibrations~$j:K→K'$ in $\ECat_\sSet$
(Bergner [\SCat, Propositions 2.3 and~2.5])
shows that their domains and codomains are cofibrant.
In the commutative square
$$\sqcd{
K&\mapright{j}&K'\cr
\mapdown{}&&\mapdown{}\cr
\UnderM\FreeM K&\mapright{\UnderM\FreeM j}&\UnderM\FreeM K'\cr
}$$
the vertical maps are unit maps for cofibrant objects $K$ and $K'$, which are weak equivalences by ^!{unit map is a weak equivalence}.
Since the map~$j$ is a weak equivalence, by the 2-out-of-3 property so is $\UnderM(\FreeM(j))$.
Since the functor~$\UnderM$ reflects weak equivalences, the map $\FreeM(j)$ is a weak equivalence.
Thus, the functor~$\FreeM$ preserves cofibrations and sends generating acyclic cofibrations to weak equivalences.
Therefore, the adjunction $\FreeM⊣\UnderM$ is a Quillen adjunction.
\ppar
To show that the adjunction $\FreeM⊣\UnderM$ is a Quillen equivalence,
consider the derived unit map of a cofibrant object $A∈\ECat_\sSet$.
Since the functor~$\UnderM$ preserves weak equivalences, the derived unit map of~$A$ can be computed as the unit map of~$A$,
which is a weak equivalence by ^!{unit map is a weak equivalence}.
\ppar
Consider the derived counit map $$ε_B:\FreeM(Q(\UnderM(B)))→B$$ of a fibrant object $B∈\ECat_\Cat$,
where $Q:\ECat_\sSet→\ECat_\sSet$ is a cofibrant replacement functor.
Since $\UnderM$ reflects weak equivalences, to show that $ε_B$ is a weak equivalence, it suffices to show that $\UnderM(ε_B)$ is a weak equivalence.
The composition of the unit map $$Q(\UnderM(B))→\UnderM(\FreeM(Q(\UnderM(B))))$$ with the map $$\UnderM(ε_B):\UnderM(\FreeM(Q(\UnderM(B))))→\UnderM(B)$$ is the cofibrant replacement map $Q(\UnderM(B))→\UnderM(B)$ of~$\UnderM(B)$,
which is a weak equivalence.
Since the unit map of the cofibrant object $Q(\UnderM(B))$ is a weak equivalence by ^!{unit map is a weak equivalence}, by the 2-out-of-3 property the map $\UnderM(ε_B)$ is a weak equivalence,
hence so is the map~$ε_B$.
\ppar
Thus, the derived unit and derived counit of the Quillen adjunction $\FreeM⊣\UnderM$ are natural weak equivalences,
hence $\FreeM⊣\UnderM$ is a Quillen equivalence.
\ppar
The case of monoids is treated in the same way (an enriched category with one object can be identified with a monoid),
using ^!{unit map is a weak equivalence} and the set of generating cofibrations described in ^!{model structure on monoids}.

\proclaim Proposition.
^^={unit map is a weak equivalence}
The unit map $X→\UnderM\FreeM X$ of a cofibrant object $X∈\ECat_\sSet$
for the adjunction $$\FreeM⊣\UnderM: \ECat_\sSet ⇄ \ECat_\Cat$$
(^!{enriched Thomason Quillen equivalence})
is a weak equivalence in $\ECat_\sSet$.
The same is true for unit maps in the adjunction $\FreeM⊣\UnderM: \Mon_\sSet ⇄ \Mon_\Cat$ (^!{monoid conventions}).

\proof Proof.
We prove by induction on a cofibrant object~$X∈\ECat_\sSet$ that the unit map $X→\UnderM\FreeM X$ of~$X$ is a weak equivalence.
If $X=∅$, the unit map of~$X$ is the identity map $∅→∅$, which is a weak equivalence.
Suppose the unit map $X→\UnderM\FreeM X$ of~$X$ is a weak equivalence
and the map $X→P$ is a cobase change of a generating cofibration~$i$ in $\ECat_\sSet$.
If $i$ is the map $∅→\{0\}$, then $P=X⊔\{0\}$ and the unit map of~$P$ is simply
the coproduct of the unit map of~$X$ and the identity map on $\{0\}$, hence is a weak equivalence.
\ppar
If $i$ is the map $T_{0,1}(δ_n)$, where $δ_n:∂Δ^n→Δ^n$ is a generating cofibration of simplicial sets,
then the cobase change square for $T_{0,1}(δ_n)$ in the category $\ECat_\sSet$
is mapped by the functor~$\FreeM$ to the cobase change square for $T_{0,1}(\fc\Sd^2(δ_n))$ in the category $\ECat_\Cat$,
since $\FreeM(T_{0,1}(δ_n))$ can be computed as $T_{0,1}(\fc\Sd^2(δ_n))$.
\ppar
Consider the unit natural transformation of commutative squares
$$\cd{
T_{0,1}(∂Δ^n)&\mapright{}&X\cr
\lmapdown{T_{0,1}(δ_n)}&&\mapdown{}\cr
T_{0,1}(Δ^n)&\mapright{}&P\cr
}
\quad⟹\qquad\quad
\cd{
\UnderM T_{0,1}(\fc\Sd^2 ∂Δ^n)&\mapright{}&\UnderM \FreeM X\cr
\lmapdown{\UnderM T_{0,1}(\fc\Sd^2 δ_n)}&&\mapdown{}\cr
\UnderM T_{0,1}(\fc\Sd^2 Δ^n)&\mapright{}&\UnderM \FreeM P\cr
}
$$
The left square is homotopy cocartesian because $T_{0,1}(δ_n)$ is a cofibration in $\ECat_\sSet$, and the latter category is left proper (^!{simplicial categories are left proper}).
The components $$T_{0,1}(∂Δ^n)→\UnderM T_{0,1}(\fc\Sd^2 ∂Δ^n),\qquad T_{0,1}(Δ^n)→T_{0,1}(\fc\Sd^2 Δ^n)$$
are weak equivalences
by ^!{properties of the Thomason adjunction},
since the functors $\UnderM$ and $\FreeM$ commute past $T_{0,1}$ in the obvious way.
The component $X→\UnderM \FreeM X$ is a weak equivalence by assumption.
Thus, the component $P→\UnderM\FreeM P$ is a weak equivalence if and only if the right square is homotopy cocartesian.
\ppar
Since the map $T_{0,1}(\fc\Sd^2 δ_n)$ is an ^{h-cofibration} in $\ECat_\Cat$ by ^!{enriched Thomason model structure},
the cobase change square
$$\cd{
T_{0,1}(\fc\Sd^2 ∂Δ^n)&\mapright{}&\FreeM X\cr
\lmapdown{T_{0,1}(\fc\Sd^2 δ_n)}&&\mapdown{}\cr
T_{0,1}(\fc\Sd^2 Δ^n)&\mapright{}&\FreeM P\cr
}$$
is homotopy cocartesian.
Thus, we have to show that the functor~$\UnderM$ sends this square to a homotopy cocartesian square.
Since the map $\UnderM T_{0,1}(\fc\Sd^2 δ_n)$ is a cofibration (hence an ^{h-cofibration}) in $\ECat_\sSet$,
we have to show that the canonical map
$$\hat P=\UnderM T_{0,1}(\fc\Sd^2 Δ^n)⊔_{\UnderM T_{0,1}(\fc\Sd^2 ∂Δ^n)}\UnderM \FreeM X→\UnderM(T_{0,1}(\fc\Sd^2 Δ^n)⊔_{T_{0,1}(\fc\Sd^2 ∂Δ^n)}\FreeM X)≅\UnderM \FreeM P$$
is a weak equivalence.
The domain and codomain of this map, when considered as objects in the undercategory of $\UnderM\FreeM X$,
admit filtrations on their hom-objects as in ^!{enriched Thomason model structure},
using the cobase changes squares for the maps $\UnderM T_{0,1}(\fc\Sd^2 δ_n)$ and $T_{0,1}(\fc\Sd^2 δ_n)$ given above,
and applying the functor~$\Ex^2\N$ to the latter filtration.
\ppar
Recall that the construction of the filtration in ^!{enriched Thomason model structure} starts by
observing that the inclusion $\Sd^2 ∂Δ^n→\Sd^2 Δ^n$ is the nerve of a map of posets $ι:A→B$,
which factors as $$\fc\Sd^2 ∂Δ^n=A\lto5{κ}W\lto5{μ}B=\fc\Sd^2 Δ^n,$$ where $W$ is the cosieve generated by~$A$ in~$B$.
We also have the induced commutative diagrams
$$
\cd{
\UnderM T_{0,1}(A)&\mapright{}&\UnderM \FreeM X\cr
\lmapdown{\UnderM T_{0,1}(κ)}&&\mapdown{}\cr
\UnderM T_{0,1}(W)&\mapright{}&\hat R\cr
\lmapdown{\UnderM T_{0,1}(μ)}&&\mapdown{}\cr
\UnderM T_{0,1}(B)&\mapright{}&\hat P,\cr
}
\qquad\qquad\qquad
\cd{
T_{0,1}(A)&\mapright{}&\FreeM X\cr
\lmapdown{T_{0,1}(κ)}&&\mapdown{}\cr
T_{0,1}(W)&\mapright{}&\bar R\cr
\lmapdown{T_{0,1}(μ)}&&\mapdown{}\cr
T_{0,1}(B)&\mapright{}&\bar P\rlap{${}=\FreeM P$,}\cr
}
\qquad\qquad\qquad
\cd{
\UnderM T_{0,1}(A)&\mapright{}&\UnderM\FreeM X\cr
\lmapdown{\UnderM T_{0,1}(κ)}&&\mapdown{}\cr
\UnderM T_{0,1}(W)&\mapright{}&\UnderM\bar R\cr
\lmapdown{\UnderM T_{0,1}(μ)}&&\mapdown{}\cr
\UnderM T_{0,1}(B)&\mapright{}&\UnderM\bar P\rlap{${}=\UnderM\FreeM P$,}\cr
}
$$
where the two pairs of squares on the left are cobase change squares and the right diagram is given by applying the functor~$\UnderM$ to the middle diagram.
\ppar
By Bergner [\SCat, Proposition~3.2], the map $\UnderM(T_{0,1}(κ))$ is an acyclic cofibration in $\ECat_\sSet$, hence the map $\UnderM\FreeM X→\hat R$ is a weak equivalence.
By ^!{enriched Thomason model structure}, the map $T_{0,1}(κ)$ is an acyclic cofibration in $\ECat_\Cat$, so the map $\FreeM X→\bar R$ is a weak equivalence.
Since the functor~$\UnderM$ preserves weak equivalences, the map $\UnderM\FreeM X→\UnderM\bar R$ is a weak equivalence in $\ECat_\sSet$.
Thus, the top squares in all three diagrams are cobase changes squares and homotopy cobase change squares.
The natural transformation from the first diagram to the third diagram has identities as its components for the left column,
as well as the top right corner.
Thus, the component $\hat R→\UnderM\bar R$ is also a weak equivalence.
\ppar
Next, consider the bottom cobase change squares in the above three diagrams.
As explained in the proof of ^!{enriched Thomason model structure},
after taking hom-objects for a fixed pair of objects $x$~and~$y$ of~$X$,
we get filtrations on
the bottom right maps $\hat R(x,y)→\hat P(x,y)$, $\bar R(x,y)→\FreeM P(x,y)$, and (after applying the functor~$\Ex^2\N$) the map $\Ex^2\N(\bar R(x,y))→\Ex^2\N(\FreeM P(x,y))$.
\ppar
For the map $\bar R(x,y)→\bar P(x,y)=\FreeM P(x,y)$ the $i$th map $\bar P_{i-1}→\bar P_i$ in the filtration is a cobase change in $\Cat$ of the map
$$β=\bar R(b,y)⨯\bar R(b,a)^{⨯(i-1)}⨯\bar R(x,a)⨯μ^{◻i}.$$
This cobase change was shown to yield a homotopy cocartesian square in ^!{enriched Thomason model structure}.
Since the functor~$\Ex^2\N$ is a right Quillen equivalence that preserves weak equivalences,
the image of the latter homotopy cocartesian square under~$\Ex^2\N$ is again homotopy cocartesian.
\ppar
Likewise, for the map $\hat R(x,y)→\hat P(x,y)$ the $i$th map $\hat P_{i-1}→\hat P_i$ in the filtration is a cobase change in $\sSet$ of the map
$$α=\hat R(b,y)⨯\hat R(b,a)^{⨯(i-1)}⨯\hat R(x,a)⨯(\Ex^2\N μ)^{◻i},$$
which is a monomorphism of simplicial sets, hence the corresponding commutative square is homotopy cocartesian.
\ppar
Thus, we have a natural transformation of homotopy cocartesian squares
$$\cd{
\dom α&\mapright{}&\hat P_{i-1}\cr
\lmapdown{α}&&\mapdown{}\cr
\codom α&\mapright{}&\hat P_i\cr
}
\quad⟹\qquad\quad
\cd{
\Ex^2\N\dom β&\mapright{}&\Ex^2\N\bar P_{i-1}\cr
\lmapdown{\Ex^2\N β}&&\mapdown{}\cr
\Ex^2\N\codom β&\mapright{}&\Ex^2\N\bar P_i.\cr
}$$
The components $\dom α→\Ex^2\N\dom β$ and $\codom α→\Ex^2\N\codom β$
are weak equivalences in $\sSet$
because the functor~$\Ex^2\N$ preserves finite products,
the map $\hat R→\UnderM\bar R$ was shown to be a weak equivalence above,
and the morphism of arrows $$(\Ex^2\N μ)^{◻i}→\Ex^2\N(μ^{◻i})$$
is a weak equivalence.
For the latter, observe that the punctured $i$-dimensional cube whose colimit yields the domain of $μ^{◻i}$
is a homotopy cocartesian cube,
and the functor $\Ex^2\N$ is a right Quillen equivalence that preserves weak equivalences,
hence it preserves homotopy colimit diagrams.
Thus, if the map $\hat P_{i-1}→\Ex^2\N\bar P_{i-1}$ is a weak equivalence,
then so is the map $\hat P_i→\Ex^2\N\bar P_i$.
\ppar
The map $\hat P_0→\Ex^2\N\bar P_0$ coincides with the map $\hat R(x,y)→\Ex^2\N(\bar R(x,y))$, which was shown to be a weak equivalence above.
By induction, the map $\hat P_i→\Ex^2\N\bar P_i$ is a weak equivalence for all $i≥0$.
Since simplicial weak equivalences are closed under filtered colimits, the map $\hat P(x,y)→\Ex^2\N(\bar P(x,y))$ is also a weak equivalence, which completes the proof.

\proclaim Remark.
The Quillen equivalence of ^!{enriched Thomason Quillen equivalence}
immediately connects (via zigzags of Quillen equivalences)
the ^{enriched Thomason model structure} on small 2-categories
to all the other models for $(∞,1)$-categories
such as quasicategories, relative categories (Barwick–Kan [\RelCat]), Segal categories, complete Segal spaces, marked simplicial sets, etc.
See Bergner [\Survey, \Book] for a review of these Quillen equivalences.

\unsection References

%$$

\yearkeytrue

\refs

\bib\EOP
Edward Szpilrajn.
Sur l'extension de l'ordre partiel.
Fundamenta Mathematicae 16 (\y{1930}), 386–389.
\doi:10.4064/fm-16-1-386-389.

\bib\AHiii[1956]
Daniel~M.~Kan.
Abstract homotopy.  III.
Proceedings of the National Academy of Sciences 42:7 (1956), 419–421.
\doi:10.1073/pnas.42.7.419.

\bib\HoAlg
Daniel~G.~Quillen.
Homotopical algebra.
Lecture Notes in Mathematics 43 (\y{1967}), Springer.
\doi:10.1007/bfb0097438.

\bib\CFHT
Peter Gabriel, Michel Zisman.
Calculus of Fractions and Homotopy Theory.
Ergebnisse der Mathematik und ihrer Grenzgebiete 35 (\y{1967}).
\doi:10.1007/978-3-642-85844-4.

\bib\EnrOp[1972]
G.~Maxwell~Kelly.
On the operads of J.~P.~May.
Reprints in Theory and Applications of Categories 13 (2005).
\http://tac.mta.ca/tac/reprints/articles/13/tr13abs.html.

\bib\HAKT
Daniel Quillen.
Higher algebraic K-theory:~I. 
In: Higher K-theories.
Lecture Notes in Mathematics 341 (\y{1973}), 85--147.
\doi:10.1007/bfb0067053.

\bib\Hocolim
Robert~W.~Thomason.
Homotopy colimits in the category of small categories.
Mathematical Proceedings of the Cambridge Philosophical Society 85:1 (\y{1979}), 91–109.
\doi:10.1017/s0305004100055535.

\bib\CatCMC
Robert~W.~Thomason.
Cat as a closed model category.
Cahiers de topologie et géométrie différentielle catégoriques 21:3 (\y{1980}), 305–324.
\numdam CTGDC_1980__21_3_305_0.

\bib\EnrCat
G.~Maxwell~Kelly.
Basic Concepts of Enriched Category Theory.
London Mathematical Society Lecture Note Series 64 (\y{1982}).
Reprints in Theory and Applications of Categories 10 (2005).
\http://tac.mta.ca/tac/reprints/articles/10/tr10abs.html.

\bib\SSCS
André Joyal, Myles Tierney.
Strong Stacks and Classifying Spaces.
In: Category Theory.
Lecture Notes in Mathematics 1488 (\y{1991}), 213–236.
\doi:10.1007/bfb0084222.

\bib\QCMSS[1993]
Sjoerd~E.~Crans.
Quillen closed model structures for sheaves.
Journal of Pure and Applied Algebra 101:1 (1995), 35–57.
\doi:10.1016/0022-4049(94)00033-f.

\bib\AMMMC[1998]
Stefan Schwede, Brooke~E.~Shipley.
Algebras and modules in monoidal model categories.
Proceedings of the London Mathematical Society 80:2 (2000), 491–511.
\arXiv:math/9801082v1, \doi:10.1112/s002461150001220x.

\bib\Replacing[1998]
Daniel Dugger.
Replacing model categories with simplicial ones.
Transactions of the American Mathematical Society 353:12 (2001), 5003–5027.
\doi:10.1090/s0002-9947-01-02661-7.

\bib\Desc[1998]
André Hirschowitz, Carlos Simpson.
Descente pour les $n$-champs.
\arXiv:math/9807049v3.

\bib\CSS[1998]
Charles Rezk.
A model for the homotopy theory of homotopy theory.
Transactions of the American Mathematical Society 353:3 (2000), 973–1007.
\arXiv:math/9811037v3, \doi:10.1090/s0002-9947-00-02653-2.

\bib\SHT
Paul~G.~Goerss, John~F.~Jardine.
Simplicial homotopy theory.
Progress in Mathematics 174 (\y{1999}).
\doi:10.1007/978-3-0346-0189-4.

\bib\MC
Mark Hovey.
Model categories.
Mathematical Surveys and Monographs 63 (\y{1999}).
\doi:10.1090/surv/063.

\bib\Dwyer
Denis-Charles Cisinski.
La classe des morphismes de Dwyer n'est pas stable par retractes.
Cahiers de topologie et géométrie différentielle catégoriques 40:3 (\y{1999}), 227–231.
\numdam CTGDC_1999__40_3_227_0.

\bib\KL
G.~Maxwell Kelly, Stephen Lack.
$\cal V$-Cat is locally presentable or locally bounded if $\cal V$ is so.
Theory and Applications of Categories 8:23 (\y{2001}), 555–575.
\eudml:122318, \http://tac.mta.ca/tac/volumes/8/n23/8-23abs.html.

\bib\EMMC[2002]
Stefan Schwede, Brooke Shipley.
Equivalences of monoidal model categories.
Algebraic \& Geometric Topology 3 (2003), 287–334.
\arXiv:math/0209342v2, \doi:10.2140/agt.2003.3.287.

\bib\MCL
Philip~S.~Hirschhorn.
Model categories and their localizations.
Mathematical Surveys and Monographs 99 (\y{2003}).
\doi:10.1090/surv/099.

\bib\Bicat[2004]
Stephen Lack.
A Quillen model structure for bicategories.
K-Theory 33:3 (2004), 185–197.
\doi:10.1007/s10977-004-6757-9.

\bib\FreePushouts[2004]
Anthony~D.~Elmendorf, Michael~A.~Mandell.
Rings, modules, and algebras in infinite loop space theory.
Advances in Mathematics 205:1 (2006), 163–228.
\arXiv:math/0403403v1, \doi:10.1016/j.aim.2005.07.007.

\bib\SCat[2004]
Julia~E.~Bergner.
A model category structure on the category of simplicial categories.
Transactions of the American Mathematical Society 359:5 (2006), 2043–2058.
\arXiv:math/0406507v2, \doi:10.1090/s0002-9947-06-03987-0.

\bib\TwoCatModel[2004]
Krzysztof Worytkiewicz, Kathryn Hess, Paul-Eugène Parent, Andrew Tonks.
A model structure a la Thomason on 2-Cat.
Journal of Pure and Applied Algebra 208:1 (2007), 205–236.
\arXiv:math/0411154v1, \doi:10.1016/j.jpaa.2005.12.010.

\bib\Three[2005]
Julia~E.~Bergner.
Three models for the homotopy theory of homotopy theories.
Topology 46 (2007), 397–436.
\arXiv:math/0504334v2, \doi:10.1016/j.top.2007.03.002.

\bib\QCSS[2006]
André Joyal, Myles Tierney.
Quasi-categories vs Segal spaces.
In: Categories in Algebra, Geometry and Mathematical Physics.
Contemporary Mathematics 431 (2007), 277–326.
\arXiv:math/0607820v2, \doi:10.1090/conm/431/08278.

\bib\Survey[2006]
Julia~E.~Bergner.
A survey of $(∞,1)$-categories.
In: Towards Higher Categories.
The IMA Volumes in Mathematics and its Applications 152 (2010), 69–83.
\arXiv:math/0610239v1, \doi:10.1007/978-1-4419-1524-5_2.

\bib\QCSC[2007]
André Joyal.
Quasi-categories vs simplicial categories.
January 7, 2007.
\https://www.math.uchicago.edu/~may/IMA/Incoming/Joyal/QvsDJan9(2007).pdf.

\bib\LR[2007]
Clark Barwick.
On left and right model categories and left and right Bousfield localizations.
Homology, Homotopy and Applications 12:2 (2010), 245–320.
\arXiv:0708.2067v2, \doi:10.4310/hha.2010.v12.n2.a9.

\bib\QCat[2008]
André Joyal.
The Theory of Quasi-Categories and its Applications.
January 31, 2008.
\https://mat.uab.cat/~kock/crm/hocat/advanced-course/Quadern45-2.pdf.

\bib\Thomncat[2008]
Thomas~M.~Fiore, Simona Paoli.
A Thomason model structure on the category of small $n$-fold categories.
Algebraic \& Geometric Topology 10:4 (2010) 1933–2008.
\arXiv:0808.4108v2, \doi:10.2140/agt.2010.10.1933.

\bib\Posets[2010]
George Raptis.
Homotopy theory of posets.
Homology, Homotopy and Applications 12:2 (2010), 211–230.
\doi:10.4310/hha.2010.v12.n2.a7.

\bib\RelCat[2010]
Clark Barwick, Daniel~M.~Kan.
Relative categories: Another model for the homotopy theory of homotopy theories.
Indagationes Mathematicae 23:1–2 (2012), 42–68.
\arXiv:1011.1691v2, \doi:10.1016/j.indag.2011.10.002.

\bib\Partial[2011]
Clark Barwick, Daniel~M.~Kan.
Partial model categories and their simplicial nerves.
\arXiv:1102.2512v2.

\bib\SOpII[2011]
Denis-Charles Cisinski, Ieke Moerdijk.
Dendroidal sets and simplicial operads.
Journal of Topology 6:3 (2013), 705–756.
\arXiv:1109.1004v1, \doi:10.1112/jtopol/jtt006.

\bib\HTEC[2012]
Fernando Muro.
Dwyer–Kan homotopy theory of enriched categories.
Journal of Topology 8:2 (2015), 377–413.
\arXiv:1201.1575v4, \doi:10.1112/jtopol/jtu029.

\bib\Lax[2012]
Jonathan Chiche.
Un théorème a de Quillen pour les 2-foncteurs lax.
Theory and Applications of Categories 30:4 (2015), 49–85.
\arXiv:1211.2319v2, \http://tac.mta.ca/tac/volumes/30/4/30-04abs.html.

\bib\HTNOii[2013]
Fernando Muro.
Homotopy theory of non-symmetric operads, II: change of base category and left properness.
Algebraic \& Geometric Topology 14 (2014), 229–281.
\arXiv:1304.6641v2, \doi:10.2140/agt.2014.14.229.

\bib\HTAPM[2013]
Michael Batanin, Clemens Berger.
Homotopy theory for algebras over polynomial monads.
Theory and Applications of Categories 32:6 (2017), 148–253.
\arXiv:1305.0086v7.

\bib\Vers[2013]
Dimitri Ara, Georges Maltsiniotis.
Vers une structure de catégorie de modèles à la Thomason sur la catégorie des $n$-catégories strictes.
Advances in Mathematics 259 (2014), 557–654.
\arXiv:1305.5086v3, \doi:10.1016/j.aim.2014.03.013.

\bib\FPMC[2014]
Lennart Meier, Viktoriya Ozornova.
Fibrancy of partial model categories.
Homology, Homotopy and Applications 17:2 (2015), 53-80.
\arXiv:1408.2743v1, \doi:10.4310/hha.2015.v17.n2.a5.

\bib\Adm[2014]
Dmitri Pavlov, Jakob Scholbach.
Admissibility and rectification of colored symmetric operads.
Journal of Topology 11:3 (2018), 559–601.
\arXiv:1410.5675v4, \doi:10.1112/topo.12008.

\bib\Properness[2014]
Philip Hackney, Marcy Robertson, Donald Yau.
Relative left properness of colored operads.
\arXiv:1411.4668v3, \doi:10.2140/agt.2016.16.2691.

\bib\FCFRC[2015]
Lennart Meier.
Fibration categories are fibrant relative categories.
Algebraic \& Geometric Topology 16:6 (2016), 3271–3300.
\arXiv:1503.02036v2, \doi:10.2140/agt.2016.16.3271.

\bib\TwoCatModelCorr[2015]
Krzysztof Worytkiewicz, Kathryn Hess, Paul-Eugène Parent, Andrew Tonks.
Corrigendum to: “A model structure à la Thomason on 2-Cat”.
Journal of Pure and Applied Algebra 220:12 (2016), 4017–4023.
\doi:10.1016/j.jpaa.2016.04.002.

\bib\Acyclic[2015]
Roman Bruckner.
A model structure on the category of small acyclic categories.
\arXiv:1508.00992v1.

\bib\Sym[2015]
Dmitri Pavlov, Jakob Scholbach.
Homotopy theory of symmetric powers.
Homology, Homotopy, and Applications 20:1 (2018), 359–397.
\arXiv:1510.04969v3, \doi:10.4310/HHA.2018.v20.n1.a20.

\bib\ThCofibrant[2016]
Roman Bruckner, Christoph Pegel.
Cofibrant objects in the Thomason model structure.
\arXiv:1603.05448v1.

\bib\AraModTh[2016]
Dimitri Ara.
Structures de catégorie de modèles à la Thomason sur la catégorie des 2-catégories strictes.
Cahiers de topologie et géométrie différentielle catégoriques 56:2 (2015), 83–108.
\arXiv:1607.03644v1, \http://cahierstgdc.com/wp-content/uploads/2017/05/DAra.pdf.

\bib\HTT
Jacob Lurie.
Higher Topos Theory.
April 9, \y{2017}.
\https://www.math.ias.edu/~lurie/papers/HTT.pdf.

\bib\Book
Julia~E.~Bergner.
The homotopy theory of $(∞,1)$-categories.
London Mathematical Society Student Texts 90 (\y{2018}), Cambridge University Press.
\doi:10.1017/9781316181874.

\bib\LAMS[2018]
Richard Garner, Magdalena Kędziorek, Emily Riehl.
Lifting accessible model structures.
Journal of Topology 13:1 (2020), 59–76.
\arXiv:1802.09889v1, \doi:10.1112/topo.12123.

\bib\HCHA
Denis-Charles Cisinski.
Higher Categories and Homotopical Algebra.
Cambridge Studies in Advanced Mathematics 180 (\y{2019}).
\doi:10.1017/9781108588737.

\bib\EICT[2022]
Emily Riehl, Dominic Verity.
Elements of ∞-Category Theory.
Cambridge Studies in Advanced Mathematics 194 (2022).
\doi:10.1017/9781108936880, \https://emilyriehl.github.io/files/elements.pdf.

%$$